\documentclass[11pt,reqno]{amsart}
\usepackage[T1]{fontenc}
\usepackage[utf8]{inputenc}
\usepackage{lmodern}
\usepackage[a4paper,margin=28mm]{geometry}
\usepackage{amsmath,amssymb,amsthm,mathtools}
\usepackage[expansion=false]{microtype}
\usepackage{enumitem}
\usepackage{tikz}
\usetikzlibrary{arrows.meta}
\usepackage{hyperref}
\hypersetup{hidelinks,pdftitle={Half-factorial generating sets and the realization of Dedekind class groups},pdfauthor={}}
\numberwithin{equation}{section}
\newtheorem{theorem}{Theorem}[section]
\newtheorem{proposition}[theorem]{Proposition}
\newtheorem{lemma}[theorem]{Lemma}
\newtheorem{corollary}[theorem]{Corollary}
\newtheorem*{maintheorem}{Main Theorem}
\theoremstyle{definition}

\newtheorem{example}[theorem]{Example}
\theoremstyle{remark}
\newtheorem{remark}[theorem]{Remark}
\newcommand{\Z}{\mathbb Z}
\newcommand{\Q}{\mathbb Q}
\newcommand{\N}{\mathbb N}
\newcommand{\Nzero}{\mathbb N_0}
\newcommand{\HF}{\mathcal H}
\newcommand{\B}{\mathcal B}
\newcommand{\F}{\mathcal F}
\newcommand{\A}{\mathcal A}
\newcommand{\LL}{\mathsf L}

\DeclareMathOperator{\ord}{ord}

\DeclareMathOperator{\q}{q}
\title[Half-factorial generating sets and class groups]{Half-factorial generating sets and the realization of Dedekind class groups}
\author[H. Kim]{Hwankoo Kim}
\address{Department of Computer Engineering, Hoseo University, Republic of Korea}
\email{hkkim@hoseo.edu}

\date{\today}

\begin{document}
\begin{abstract}
We prove a cardinality obstruction to realizing an abelian group as
the ideal class group of a half-factorial Dedekind domain. If an
abelian $p$-group $P$ satisfies $\bigcap_{r\ge1}p^rP=0$ and admits a
half-factorial generating subset, then
$|P|\le\max\{\aleph_0,|P/pP|\}$. The proof separates generators of
equal order by finite quotients. Its group-theoretic ingredient is
a direct summand generated by selected members of an arbitrary
finite generating set. The torsion subgroup of
$\prod_{j\ge1}\Z/p^j\Z$ is uncountable and $p$-adically separated,
but has countable quotient modulo $p$; hence it has no
half-factorial generating subset. This gives a negative answer to
the universal half-factorial Dedekind realization problem for every
prime. Pure subgroups give counterexamples of every uncountable
cardinality at most the continuum, so the least cardinality of a
counterexample is $\aleph_1$. We also study free-kernel lifting,
order-preserving lifts across torsion kernels, character-compatible
generators, and the failure of compatibility under unions and
contraction. These results distinguish global realizability from
half-factoriality on bounded subgroups and from integrality of
cross numbers.
\end{abstract}
\maketitle

\section{Introduction}\label{sec:introduction}

Factorization into irreducibles is unique in a factorial domain, but
uniqueness of the number of factors is a substantially weaker
condition. A domain is \emph{half-factorial} if every nonzero nonunit
factors into irreducibles and all such factorizations have the same
length. Carlitz's theorem characterizes half-factorial rings of
integers of algebraic number fields by class number at most two
\cite{Carlitz}. The corresponding question for arbitrary Dedekind
domains is different because the distribution of prime ideals among
the ideal classes is no longer fixed by the number-field setting.

Claborn proved that every abelian group is the ideal class group of
a Dedekind domain \cite{Claborn}. Grams and subsequent realization
theorems permit prescribed classes of prime divisors under the
appropriate generating condition
\cite{Grams,ChapmanGeroldinger}. Thus the abstract class group does
not by itself determine the factorization lengths of a particular
Dedekind domain. Nevertheless, it is natural to ask whether every
abelian group admits \emph{some} prime-divisor distribution for
which the resulting domain is half-factorial. Equivalently, one asks:
\begin{quote}
Is every abelian group isomorphic to the ideal class group of a
half-factorial Dedekind domain?
\end{quote}
This is Problem~15 in the survey of Geroldinger, Kim, and Loper
\cite[\S12.1]{GKL}.

The group-theoretic formulation uses zero-sum sequences. For
$S\subseteq G$, let $\B(S)$ be the monoid of finite sequences over
$S$ whose sum is zero. The classical realization criterion says
that $G$ occurs as the class group of a half-factorial Dedekind
domain precisely when it has a monoid-generating subset $S$ for
which $\B(S)$ is half-factorial
\cite[Theorem~3.6]{ChapmanGeroldinger}. Group generation and monoid
generation agree for torsion groups. For such a group, the
Skula--\'{S}liwa cross-number criterion characterizes half-factoriality
by $\mathsf k(U)=1$ for every minimal nonempty zero-sum sequence $U$,
where $\mathsf k(U)=\sum_s v_s(U)/\ord(s)$
\cite[Lemma~3.2]{GaoGeroldinger}.

Michel and Steffan obtained early positive realization results
\cite{MichelSteffan}. Geroldinger and G\"obel proved that every
Warfield group admits half-factorial generators and, in particular,
that every countable abelian group does so
\cite[Theorems~5.1 and~5.3]{GeroldingerGoebel}. They also reduced
the universal problem to reduced uncountable unbounded primary
groups and singled out the torsion-completion of
$\bigoplus_{j\ge1}C_{p^j}$ as an unresolved example
\cite[pp.~855--856]{GeroldingerGoebel}. The positive lifting and
reduction statements are established results; we use them with
explicit attribution rather than treating them as contributions
of this paper.

Our main contribution is a restriction on \emph{arbitrary}
half-factorial generating sets. It does not require generators of
the socle, compatible generators of bounded layers, or a prescribed
character. Say that an abelian $p$-group is \emph{$p$-adically
separated} if $\bigcap_{r\ge1}p^rP=0$. The principal conclusions
are as follows.

\begin{maintheorem}
Let $p$ be a prime.
\begin{enumerate}[label=\textup{(\arabic*)},leftmargin=2.2em]
\item If a $p$-adically separated abelian $p$-group $P$ has a
half-factorial generating subset, then
\[
 |P|\le\max\{\aleph_0,|P/pP|\}.
\]
\item The group $P_p=t(\prod_{j\ge1}\Z/p^j\Z)$, where $t$ denotes
the torsion subgroup, has no half-factorial generating subset.
Consequently it is not the ideal class group of a half-factorial
Dedekind domain, nor the divisor class group of a half-factorial
Krull monoid.
\item For every uncountable cardinal $\kappa\le2^{\aleph_0}$,
there is a reduced abelian $p$-group of cardinality $\kappa$ with
no half-factorial generating subset. The least cardinality of an
abelian group not realizable as the class group of a half-factorial
Dedekind domain is therefore $\aleph_1$.
\end{enumerate}
\end{maintheorem}

The proof of the cardinality bound is finite at its decisive step.
Given a finite abelian $p$-group $G$, a generating set $S$, and
$r\ge1$, we select $A\subseteq S$ such that $K=\langle A\rangle$
is a direct summand with $p^rG\subseteq K$ and
$K[p]\subseteq p^rG$. Two distinct generators of the same order
that become equal modulo $K$ without losing order yield a minimal
zero-sum sequence of cross number greater than one. It follows that,
for a half-factorial generating set of any abelian $p$-group,
\[
 \ord(s)=\ord(t)=p^n,\quad s-t\in p^rP,\quad
 p^{n-1}s\notin p^rP \quad\Longrightarrow\quad s=t.
\]
Stratifying the generators by $n$ and $r$ gives the cardinality
bound. This separates the present argument from obstructions to
particular extension procedures: the necessary subgroup is selected
from the very generating set being tested.

The complementary results address how half-factorial generating
sets behave under algebraic constructions. We extend the
Geroldinger--G\"obel free-kernel decomposition to quotients that may
contain elements of infinite order. For torsion kernels, we prove
that half-factoriality of an order-preserving lifted generating
set is equivalent to splitting along those representatives. We
also construct weakly half-factorial generators for every torsion
abelian group, impose a prescribed character in bounded primary
groups, and give explicit obstructions to thinning, adjoining
roots, and repeated multiplication by a prime. The weak character
criterion and its finite-group antecedents belong to the established
theory \cite{GaoGeroldinger,RadziejewskiSchmid,Schmid2005}; the
additional constructions and proofs are presented separately from
the main nonrealization argument.

Section~\ref{sec:preliminaries} fixes conventions and recalls the
realization and cross-number criteria. Section~\ref{sec:cardinality}
proves the Main Theorem and its cardinal refinements.
Section~\ref{sec:extensions} treats lifting, the classical reduction,
and pure extensions. Section~\ref{sec:characters} studies characters
and finite obstructions. Section~\ref{sec:socle} concludes with
socle compatibility, contraction, and the limitations of the
cardinality criterion.

\section{Zero-sum sequences and realization}\label{sec:preliminaries}

This section establishes the conventions used throughout the paper,
records the exact connection with Dedekind class groups, and gives
the elementary positive examples needed later. We distinguish the
ambient group of labels from the intrinsic divisor class group of
an unlabelled block monoid. Standard references for the group and
factorization theory are \cite{Fuchs,FuchsII,GHK}.

\subsection{Conventions and factorization lengths}

We write $\N=\{1,2,\ldots\}$ and $\Nzero=\N\cup\{0\}$.
All groups are abelian and written additively. All monoids are
commutative and cancellative and have an identity. A monoid is
\emph{reduced} if its only unit is its identity; an \emph{atom} is
a nonunit that is not a product of two nonunits. A monoid is
\emph{atomic} if every nonunit is a finite product of atoms.
Write $\A(M)$ for the set of atoms of a monoid $M$ and $\q(M)$
for its quotient group. For an atomic monoid $M$ and $a\in M$,
let $\LL_M(a)$ be the set
of lengths of its factorizations into atoms, and set
$\LL_M(a)=\{0\}$ for a unit. The monoid is \emph{half-factorial}
if every such set is a singleton. A domain has the corresponding
property when its multiplicative monoid of nonzero elements does.
All rings are commutative with $1\ne0$, and ring homomorphisms
preserve the identity.

For a group $G$ and $S\subseteq G$, write $\langle S\rangle$ for
the generated subgroup and $[S]$ for the generated submonoid,
including the empty sum. Thus $[S]=G$ requires nonnegative
coefficients. If $G$ is torsion, then $[S]=\langle S\rangle$.
Write $C_n=\Z/n\Z$ and identify $C_{p^\infty}$ with the
$p$-primary subgroup of $\Q/\Z$. An abelian $p$-group is a group
all of whose elements have finite $p$-power order. For $m\ge1$,
put $G[m]=\{g\in G:mg=0\}$. The subgroup $P[p]$ of a $p$-group
is its \emph{socle}, an $\mathbb F_p$-vector space.
A group is \emph{divisible} if multiplication by every positive
integer is surjective, and \emph{reduced} if it contains no nonzero
divisible subgroup. A subgroup $H\subseteq G$ is \emph{pure} if
$H\cap mG=mH$ for every $m\ge1$. A $p$-group is \emph{bounded}
if $p^NP=0$ for some $N$, and \emph{unbounded} otherwise.
The axiom of choice is used for bases and cardinal comparisons.

Let $\F(S)$ be the free commutative monoid on $S$. Its members
are finite sequences, written $B=\prod_{s\in S}s^{v_s(B)}$ with
finite support. The sequence length and sum are
$|B|=\sum_s v_s(B)$ and $\sigma(B)=\sum_s v_s(B)s$.
A \emph{subsequence} means a divisor in $\F(S)$, including multiplicities. A \emph{zero-sum sequence} over an additive abelian group $G$ is a finite sequence of elements of $G$, with repetitions allowed and order disregarded, whose sum is $0$. It is \emph{minimal} if it is nonempty and no proper nonempty subsequence has sum $0$.
The \emph{block monoid} is
\[
 \B(S)=\{B\in\F(S):\sigma(B)=0\}.
\]
Its identity is the empty sequence. Its atoms are exactly the
minimal nonempty zero-sum sequences; induction on $|B|$ proves
atomicity. A subset $S$ is \emph{half-factorial} if $\B(S)$ is
half-factorial. We write $\HF$ for the class of groups possessing
a half-factorial subset that generates them as a monoid.
The empty set is allowed for the trivial group.

Adding or deleting the letter zero does not affect generation or
half-factoriality: if $0\notin S$, then
$\B(S\cup\{0\})\cong\B(S)\times\Nzero$.
Unless stated otherwise, generating sets omit zero. We use the
usual convention $\ord(0)=1$ when zero occurs in a sequence.
Half-factoriality has finite character: a set is half-factorial
precisely when each of its finite subsets is half-factorial
\cite[Lemma~3.1]{GaoGeroldinger}. In particular, subsets and unions
of nested chains of half-factorial sets are half-factorial.
Indeed, any failure is witnessed by two factorizations using only
finitely many letters. This observation concerns chains of sets,
not arbitrary increasing unions of the groups they generate.

\subsection{The realization and cross-number criteria}

A Dedekind domain is a Noetherian integrally closed domain of
Krull dimension at most one. Its ideal class group is the group
of nonzero fractional ideals modulo nonzero principal fractional
ideals. 
Let $H$ be a commutative cancellative monoid and let
$F=\mathcal{F}(P)$ be the free commutative monoid on a set $P$.
A monoid homomorphism $\varphi:H\to F$ is a
\emph{divisor homomorphism} if, for all $a,b\in H$,
$\varphi(a)\mid_F\varphi(b)$ implies $a\mid_H b$.
It is a \emph{divisor theory} if, in addition, for every
$f\in F$ there exist $m\ge1$ and $a_1,\ldots,a_m\in H$ such that
\[
 f=\gcd_F\bigl(\varphi(a_1),\ldots,\varphi(a_m)\bigr).
\]
Here greatest common divisors in $F$ are computed by taking
the minimum exponent of each generator in $P$.
A commutative cancellative monoid is a \emph{Krull monoid}
if it admits a divisor theory. For a divisor theory
$\varphi:H\to F$, its \emph{divisor class group} is
\[
 \mathcal{C}(H)
 =\operatorname{q}(F)\big/
  \operatorname{q}\bigl(\varphi(H)\bigr).
\]
This group is independent, up to
isomorphism, of the chosen divisor theory; see \cite{GHK}.
We recall the realization theorem
in the form needed here. The existence statement is classical
\cite{Grams,ChapmanGeroldinger}; its formulation in terms of
generating subsets is also recorded in
\cite[p.~598]{GaoGeroldinger} and
\cite[p.~842]{GeroldingerGoebel}.

\begin{theorem}[Classical realization criterion]\label{thm:realization}
Let $G$ be an abelian group. The following are equivalent:
\begin{enumerate}[label=\textup{(\arabic*)},leftmargin=2.2em]
\item $G$ has a half-factorial subset $S$ with $[S]=G$;
\item $G$ is isomorphic to the ideal class group of a half-factorial
Dedekind domain;
\item $G$ is isomorphic to the divisor class group of a half-factorial
Krull monoid.
\end{enumerate}
More precisely, for every nonempty $S\subseteq G$ with $[S]=G$,
there is a Dedekind domain whose class group is identified with $G$
and whose classes containing nonzero prime ideals are exactly $S$.
Its factorization lengths are transferred to $\B(S)$, so it is
half-factorial precisely when $S$ is half-factorial.
\end{theorem}

\begin{proof}[Explanation of the equivalence]
The prescribed prime-class realization theorem supplies the domain.
If $(a)$ factors into prime ideals, replace each prime ideal by its
class to obtain a block over $S$. A zero-sum subsequence corresponds
to a principal subproduct; its generator gives the corresponding
element factor. Partitioning the prime-ideal occurrences lifts
factorizations, including their lengths. Conversely, irreducible
elements give minimal zero-sum subsequences under this transfer.
These are the standard block-transfer properties
\cite{ChapmanGeroldinger,GHK}.

The same argument applies to a divisor theory of a Krull monoid.
The classes of its prime divisors generate the class group as a
monoid, and the block transfer preserves sets of lengths. Thus (3)
implies (1), while (1) implies (2) by realization and (2) implies
(3) by the divisor theory of a Dedekind domain. The trivial group
may be realized by $k[X]$ for any field $k$.
\end{proof}

For a subset $T$ of a torsion group and $B\in\F(T)$, define its
\emph{cross number} by
\[
 \mathsf k(B)=\sum_{t\in T}\frac{v_t(B)}{\ord(t)}.
\]
The following criterion is due to Skula and \'{S}liwa; see
\cite[Lemma~3.2]{GaoGeroldinger} and
\cite[Theorem~4.1]{Schmid2005}. We include its short proof because
it is the only numerical factorization criterion used in the
nonrealization argument.
\begin{proposition}[Cross-number criterion]\label{prop:cross}
Let $G$ be a torsion abelian group and let $T\subseteq G$, with
$\ord(0)=1$ if $0\in T$. Then $T$ is
half-factorial if and only if $\mathsf k(U)=1$ for every
$U\in\A(\B(T))$.
\end{proposition}

\begin{proof}
If every atom has cross number one, additivity of $\mathsf k$
shows that all factorizations of a fixed block have the same
length. Conversely, assume half-factoriality and let $U$ be an
atom. Choose a positive integer $m$ divisible by the orders of
all elements in its support. The block $U^m$ has a factorization
of length $m$ using $U$, and a factorization using the pure atoms
$t^{\ord(t)}$ of length
$\sum_t mv_t(U)/\ord(t)=m\mathsf k(U)$. Equality of these
lengths gives $\mathsf k(U)=1$.
\end{proof}

\subsection{Labels and the intrinsic class group}

The quotient of the ambient free group by the group of blocks need
not be the intrinsic divisor class group of an unlabelled block
monoid. The distinction is relevant when applying a realization theorem;
compare the divisor-theoretic discussion in
\cite[\S2]{Schmid2010}.

\begin{example}\label{ex:class-trap}
Let $G=C_n=\langle g\rangle$, with $n\ge2$, and let $S=\{g\}$.
Then $[S]=G$ and $\B(S)=\{g^{nk}:k\in\Nzero\}$ is a free monoid
on one generator. Its intrinsic divisor class group is therefore
zero. On the other hand,
$\q(\F(S))/\q(\B(S))\cong\Z/n\Z$. The inclusion
$\B(S)\subseteq\F(S)$ is saturated, but is not a divisor theory.
Indeed, $g$ is not a greatest common divisor of finitely many
members of $\B(S)$.
\end{example}

An inclusion $H\subseteq\mathcal{F}(P)$ is \emph{saturated} if
$H=\operatorname{q}(H)\cap\mathcal{F}(P)$ inside
$\operatorname{q}(\mathcal{F}(P))$; equivalently, the inclusion
is a divisor homomorphism. For such an inclusion to be a
divisor theory, it suffices that every $p\in P$ be the
greatest common divisor of finitely many elements of $H$.
For such an inclusion the divisor class group is
$\q(\F(P))/\q(H)$.

\begin{proposition}\label{prop:labels}
Let $G$ be an abelian group and let $S\subseteq G$ satisfy $[S]=G$.
Set $P=S\times\N$, let $\beta:\F(P)\to G$ send $(s,j)$ to $s$,
and put $H=\beta^{-1}(0)$. Then $H\subseteq\F(P)$ is a divisor
theory with class group $G$. Forgetting the second coordinate gives
a surjective homomorphism $H\to\B(S)$ that preserves all sets of
factorization lengths. In particular, $H$ is half-factorial if and
only if $S$ is half-factorial.
\end{proposition}

\begin{proof}
The inclusion is saturated because the quotient of two blocks, when
it is a sequence, still has sum zero. Fix a label $p=(s,j)$. Since
$[S]=G$, choose a sequence with sum $-s$. Give its letters labels in
two disjoint sets, both avoiding $p$, to obtain $W_1,W_2$. Then
$pW_1,pW_2\in H$ and their greatest common divisor in $\F(P)$ is
$p$. This also works for $s=0$, using the empty sequence. Thus the
inclusion is a divisor theory.

Extend $\beta$ to $\q(\F(P))$. To verify its kernel explicitly,
let $z\in\ker\beta$. Choose a sequence $M$ clearing all negative
coordinates of $z$, so that $zM\in\F(P)$. Choose a further
sequence $C$ with $\beta(C)=-\beta(M)$. Then $MC,zMC\in H$ and
$z=(zMC)/(MC)\in\q(H)$. Since $\beta$ is surjective, the quotient
group is $G$.

Every block over $S$ can be labelled, so the forgetting map is
surjective. A factorization of its image can be lifted by
partitioning the occurrences of each class among the factors,
retaining their original labels. Conversely, forgetting labels
sends atoms to atoms: a proper zero-sum subsequence of the image
would lift to a proper zero-sum subsequence of the original block.
The map therefore preserves exactly the factorization lengths.
\end{proof}

The labelled construction identifies the divisor class group and the
length-preserving map explicitly. Its Dedekind-domain interpretation
uses Theorem~\ref{thm:realization}; no identification of an arbitrary
Krull monoid with a domain's reduced multiplicative monoid is needed.

\subsection{Elementary positive constructions}

Direct sums provide the elementary positive examples used later.
The block decomposition is standard
\cite[Lemma~3.5]{GeroldingerGoebel}; it is a restricted product,
so each sequence involves only finitely many components.

\begin{lemma}\label{lem:directsums}
Let $(G_j)_{j\in J}$ be abelian groups and suppose that
$T_j\subseteq G_j\setminus\{0\}$ is half-factorial with
$[T_j]=G_j$. Embed the $T_j$ in $G=\bigoplus_{j\in J}G_j$ and
put $T=\bigcup_j T_j$. Then $[T]=G$ and
$\B(T)\cong\bigoplus_{j\in J}\B(T_j)$, where the right side is
the restricted direct product of monoids. In particular,
$\HF$ is closed under arbitrary direct sums.
\end{lemma}

\begin{proof}
Both group elements and sequences have finite support. A sequence
on $T$ has sum zero precisely when its component in each $T_j$
has sum zero. Consequently each atom belongs to a single
component, and every factorization length is the sum of the
uniquely determined component lengths.
\end{proof}

The elementary cyclic examples already give all direct sums of
cyclic groups. For a finite cyclic group generated by $c$ of order
$n\ge2$, the set $\{c\}$ has the single block atom $c^n$.
For the infinite cyclic group generated by $e$, the set $\{e,-e\}$
has the single block atom $e(-e)$. These familiar examples and
Lemma~\ref{lem:directsums} show that all direct sums of cyclic
groups lie in $\HF$; see also \cite[p.~598]{GaoGeroldinger}.

The divisible primary example is equally explicit. The existence
assertion below is classical; compare
\cite[Corollary~5.4]{GaoGeroldinger}. We include the short argument
because it supplies the divisible part in the reduction theorem.

\begin{proposition}\label{prop:prufer}
For every prime $p$, the Pr\"ufer group $C_{p^\infty}$ belongs to
$\HF$. More precisely, identifying it with the $p$-primary subgroup
of $\Q/\Z$, the set
$T=\{p^{-n}+\Z:n\ge1\}$ is monoid-generating and half-factorial.
\end{proposition}

\begin{proof}
Put $c_n=p^{-n}+\Z$. The $c_n$ generate the group; because it is
torsion, they also generate it as a monoid. Give $c_n$ weight
$p^{-n}$. A sequence is a block exactly when its total weight is
an integer.

A finite collection of these weights having positive integer sum
$m$ can be partitioned into $m$ collections of weight one. Indeed,
let $N$ be the largest occurring denominator exponent. For $N\ge1$,
the number of occurrences of $p^{-N}$ is divisible by $p$: multiply
the total weight by $p^N$ and reduce modulo $p$. Bundle these
occurrences into packets of $p$ and regard each packet as a weight
$p^{-(N-1)}$. Repeat. This process ends with $m$ packets of weight
one, each recording a disjoint subsequence of the original sequence.

It follows that every block of weight greater than one is
reducible. Conversely, a block of weight one has no proper
nonempty block subsequence, because any such subsequence has
positive integral weight. Thus the atoms have weight exactly one,
and the total weight of a block is the length of each of its
factorizations.
\end{proof}

We will also use the established countable case:
every countable abelian group belongs to $\HF$
\cite[Theorem~5.3, with $F=0$]{GeroldingerGoebel}.
Every bounded torsion abelian group belongs to $\HF$ as well,
because it is a direct sum of finite cyclic groups \cite{Fuchs}.
These positive cases concern existence of a suitable generating
set, not the extendibility of an arbitrarily prescribed one.

\section{Separation and a cardinality obstruction}\label{sec:cardinality}

We first construct a finite zero-sum obstruction from coincident
quotient images. We then select a suitable direct summand from an
arbitrary finite generating set and derive the separation theorem.
The resulting cardinality bound excludes the torsion-completion and
its uncountable pure subgroups. The proof uses neither the later
lifting results nor a hypothesis about generators of the socle.

\subsection{Collisions and a finite separation subgroup}

The first observation turns two order-preserving representatives of
one quotient element into an atom of cross number greater than one.
No half-factoriality hypothesis is needed to construct this atom.

\begin{lemma}\label{lem:collision}
Let $G$ be a torsion abelian group, let
$S\subseteq G\setminus\{0\}$, and let $A\subseteq S$. Put
$K=\langle A\rangle$. Suppose that distinct $s,t\in S\setminus K$
satisfy
\begin{equation}\label{eq:collision-hypothesis}
 s+K=t+K,\qquad
 \ord_G(s)=\ord_G(t)=\ord_{G/K}(s+K)=q.
\end{equation}
Then there is a nonempty sequence $C\in\F(A)$ such that
$U=st^{q-1}C$ is an atom of $\B(S)$. In particular,
$\mathsf k(U)=1+\mathsf k(C)>1$, so $S$ is not half-factorial.
\end{lemma}

\begin{proof}
Since $G$ is torsion, $A$ generates $K$ as a monoid. Choose a
sequence $C$ over $A$ of least possible sequence length subject to
$\sigma(C)=t-s$. It is nonempty because $s\ne t$. The sequence
$U=st^{q-1}C$ has sum zero.

Let $V$ be a nonempty zero-sum subsequence of $U$. The letters from
$C$ vanish in $G/K$, while all $q$ occurrences in $st^{q-1}$ have
the same image of order $q$. Thus $V$ contains either none or all
of these $q$ occurrences. In the first case $V$ is a nonempty block
subsequence of $C$; deleting it gives a shorter sequence with sum
$t-s$, contrary to the choice of $C$. In the second case
$V=st^{q-1}C'$ for a subsequence $C'$ of $C$, and
$\sigma(C')=t-s$. Minimality forces $C'=C$. Therefore $U$ is an
atom. Its first $q$ letters have cross number one by
\eqref{eq:collision-hypothesis}, and $C$ has positive cross number.
Apply Proposition~\ref{prop:cross}.
\end{proof}

The next lemma is purely group-theoretic. Its generating subset is
selected from the given set, a point essential to the subsequent
application of Lemma~\ref{lem:collision}.

\begin{lemma}\label{lem:finite-separator}
Let $p$ be a prime, let $G$ be a finite abelian $p$-group, let
$S\subseteq G$ generate $G$, and let $r\ge1$. There is a subset
$A\subseteq S$ such that, for $K=\langle A\rangle$,
\begin{equation}\label{eq:finite-separator}
 p^rG\subseteq K,\qquad K[p]\subseteq p^rG.
\end{equation}
In fact, $K$ can be chosen as a direct summand isomorphic to the
sum of the cyclic factors of $G$ whose orders exceed $p^r$.
Here $K[p]=\{x\in K:px=0\}$.
\end{lemma}

\begin{proof}
\emph{Step 1: separate the cyclic factors by exponent.}
By the finite abelian-group structure theorem \cite{Fuchs}, write
\[
 G=L\oplus H,\qquad p^rL=0,\qquad
 H=\bigoplus_{j=1}^{d}\langle h_j^0\rangle,
 \quad \ord(h_j^0)=p^{a_j},\quad a_j\ge r+1.
\]
If $d=0$, take $A=\varnothing$ and $K=0$. Suppose $d>0$.
The projections of $S$ generate $H$, so choose distinct
$u_1,\ldots,u_d\in S$ whose projections form a basis of $H/pH$.
Write $u_i=\ell_i+h_i$ with $\ell_i\in L$ and $h_i\in H$.

\emph{Step 2: show that the selected subgroup is a graph.}
The elements $h_i$ generate $H$. Indeed, for
$H_0=\langle h_1,\ldots,h_d\rangle$ one has $H=H_0+pH$;
iteration and the finite exponent of $H$ give $H=H_0$.
Moreover,
\[
 H/p^{r+1}H\cong(\Z/p^{r+1}\Z)^d.
\]
The images of the $h_i$ constitute a basis of this module: they
generate it, and the resulting surjection from
$(\Z/p^{r+1}\Z)^d$ is an isomorphism, since source and codomain
have the same finite cardinality. Consequently, every integer
relation $\sum_i c_ih_i=0$ satisfies
$p^{r+1}\mid c_i$ for every $i$. Since $p^rL=0$, the same
coefficients give $\sum_i c_i\ell_i=0$.

It follows that the assignments $h_i\mapsto\ell_i$ define a
homomorphism $\psi:H\to L$. The subgroup
$K=\langle u_1,\ldots,u_d\rangle$ is its graph:
\[
 K=\{\psi(h)+h:h\in H\},\qquad G=L\oplus K.
\]
\emph{Step 3: verify the two power conditions.}
For every $h\in H$, $p^r\psi(h)=0$, so
$p^rK=p^rH=p^rG$, proving the first inclusion in
\eqref{eq:finite-separator}. Also $H[p]\subseteq p^rH$,
because all cyclic factors of $H$ have order at least $p^{r+1}$.
The homomorphism $\psi$ vanishes on $p^rH$, and hence on $H[p]$.
Thus $K[p]=H[p]\subseteq p^rG$, with $H$ regarded as the chosen
summand of $G$. Take $A=\{u_1,\ldots,u_d\}$.
\end{proof}

\subsection{Separation of arbitrary half-factorial generators}

Only finitely many generators are needed to witness a relation
$s-t\in p^rP$. This permits the finite lemma to be used in groups
of arbitrary cardinality, without a pure-subgroup assumption.

\begin{theorem}\label{thm:generator-separation}
Let $p$ be a prime, let $P$ be an abelian $p$-group, and let
$S\subseteq P\setminus\{0\}$ be a half-factorial generating set.
Suppose that $s,t\in S$ and $n,r\ge1$ satisfy
\begin{equation}\label{eq:generator-separation}
 \ord(s)=\ord(t)=p^n,\qquad
 s-t\in p^rP,\qquad p^{n-1}s\notin p^rP.
\end{equation}
Then $s=t$.
\end{theorem}

\begin{proof}
\emph{Step 1: choose a finite witness for divisibility.}
Suppose that $s\ne t$. Choose $y\in P$ such
that $s-t=p^ry$. Since $S$ generates $P$, there is a finite
subset $S_0\subseteq S$ containing $s,t$ with
$y\in G_0=\langle S_0\rangle$. The group $G_0$ is a finite
abelian $p$-group. Orders of its elements are the same as in $P$,
and
\[
 s-t\in p^rG_0,\qquad p^{n-1}s\notin p^rG_0.
\]
\emph{Step 2: preserve the orders in a quotient.}
Apply Lemma~\ref{lem:finite-separator} to $G_0,S_0,r$, obtaining
$A\subseteq S_0$ and $K=\langle A\rangle$ with
$p^rG_0\subseteq K$ and $K[p]\subseteq p^rG_0$.

The element $p^{n-1}s$ has order $p$ and does not belong to $K$;
otherwise it would lie in $K[p]\subseteq p^rG_0$. Every nonzero
subgroup of the cyclic group $\langle s\rangle$ contains
$p^{n-1}s$. Therefore $\langle s\rangle\cap K=0$, and $s+K$
has order $p^n$ in $G_0/K$. Since $s-t\in p^rG_0\subseteq K$,
the elements $s,t$ have the same nonzero image in this quotient.

\emph{Step 3: produce a large-cross-number atom.}
Lemma~\ref{lem:collision} gives an atom over $S_0$ of cross number
greater than one. It remains an atom over $S$, because its
subsequences have the same sums in $G_0$ and $P$.
This contradicts Proposition~\ref{prop:cross}.
\end{proof}

For an abelian $p$-group $P$, the \emph{$p$-adic topology} is the
group topology having $\{p^rP:r\ge1\}$ as a neighbourhood basis
at zero; see \cite[Example~4.1(a)]{DikranjanGiordanoBruno}.
We call $P$ \emph{$p$-adically separated} if
$\bigcap_{r\ge1}p^rP=0$, or equivalently, if its $p$-adic
topology is Hausdorff. The separation theorem bounds the number of generators at each
order and divisibility level. Here and below, cardinal inequalities
are taken in the usual set-theoretic sense.

\begin{theorem}\label{thm:cardinality}
Let $p$ be a prime and let $P$ be a $p$-adically separated abelian
$p$-group. If $P$ has a half-factorial generating subset, then
\begin{equation}\label{eq:cardinality-bound}
 |P|\le\max\{\aleph_0,|P/pP|\}.
\end{equation}
In particular, a $p$-adically separated abelian $p$-group with
countable $P/pP$ can admit a half-factorial generating subset only
if $P$ is countable.
\end{theorem}

\begin{proof}
Let $S\subseteq P\setminus\{0\}$ be a half-factorial generating
set. For positive integers $n,r$, put
\[
 S_{n,r}=\{s\in S:\ord(s)=p^n\text{ and }
                      p^{n-1}s\notin p^rP\}.
\]
The natural map $S_{n,r}\to P/p^rP$ is injective by
Theorem~\ref{thm:generator-separation}. Every $s\in S$ belongs
to some $S_{n,r}$: if $\ord(s)=p^n$, its nonzero multiple
$p^{n-1}s$ cannot belong to all $p^rP$. Thus
\begin{equation}\label{eq:generator-bound}
 S=\bigcup_{n,r\ge1}S_{n,r},\qquad
 |S_{n,r}|\le |P/p^rP|.
\end{equation}

Set $\kappa=\max\{\aleph_0,|P/pP|\}$. Multiplication by $p^j$
induces a surjection $P/pP\to p^jP/p^{j+1}P$ for every $j\ge0$.
The finite filtration of $P/p^rP$ with these successive quotients
therefore gives $|P/p^rP|\le\kappa$. Equation
\eqref{eq:generator-bound} yields $|S|\le\kappa$, and the group
generated by at most $\kappa$ elements has cardinality at most
$\kappa$. This proves \eqref{eq:cardinality-bound}.
\end{proof}

A countable dense subset meets every open coset, so the same
obstruction has a topological formulation.

\begin{corollary}\label{cor:separable-topology}
Let $p$ be a prime and let $P$ be an abelian $p$-group whose $p$-adic topology is Hausdorff
and has a countable dense subset. If $P$ admits a half-factorial
generating subset, then $P$ is countable.
\end{corollary}

\begin{proof}
Every coset of $pP$ is open and must meet a countable dense subset.
Hence $P/pP$ is countable. Apply Theorem~\ref{thm:cardinality}.
\end{proof}

\subsection{The torsion-completion counterexample}

The following group is the explicit test case proposed in
\cite[p.~855]{GeroldingerGoebel}. We verify directly all group
properties needed for the obstruction. Its identification as the
torsion-completion of $\bigoplus_{j\ge1}C_{p^j}$ is also discussed
in Proposition~\ref{prop:completion} below.

\begin{theorem}\label{thm:counterexample}
For every prime $p$, let
\begin{equation}\label{eq:counterexample-group}
 P_p=t\!\left(\prod_{j\ge1}\Z/p^j\Z\right)
 =\left\{(x_j)_j:
       p^m x_j=0\text{ for all }j\text{ and some }m\ge0\right\},
\end{equation}
where $t(-)$ denotes the torsion subgroup. Then $P_p$ is a reduced,
unbounded abelian $p$-group of cardinality $2^{\aleph_0}$ and
\begin{equation}\label{eq:completion-quotients}
 \bigcap_{r\ge1}p^rP_p=0,\qquad
 P_p/p^rP_p\cong
 \bigoplus_{j\ge1}\Z/p^{\min(j,r)}\Z
 \quad(r\ge1).
\end{equation}
It has no half-factorial generating subset. Consequently $P_p$ is
not isomorphic to the ideal class group of any half-factorial
Dedekind domain, so the universal half-factorial Dedekind realization problem
has a negative answer. The group $P_p$ is not the divisor class
group of a half-factorial Krull monoid either.
\end{theorem}

\begin{proof}
\emph{Step 1: cardinality and separation.}
Write $C_j=\Z/p^j\Z$. The description in
\eqref{eq:counterexample-group} shows that $P_p$ is a $p$-group.
It contains $\bigoplus_{j\ge1}C_j$, so its exponent is unbounded.
Its order-$p$ subgroup is
\[
 P_p[p]=\prod_{j\ge1}C_j[p]\cong\prod_{j\ge1}\Z/p\Z.
\]
This subgroup has cardinality $2^{\aleph_0}$, and the ambient
countable product of finite groups has the same cardinality.
Thus $|P_p|=2^{\aleph_0}$.

If $x\in\bigcap_{r\ge1}p^rP_p$, then for each $j$ the
$j$th coordinate of $x$ lies in $p^jC_j=0$. Hence $x=0$.
This also proves reducedness: every divisible subgroup is contained
in every $p^rP_p$.

\emph{Step 2: compute all finite power quotients.}
Fix $r\ge1$. Coordinate reduction induces a homomorphism
\[
 \rho_r:P_p\longrightarrow\bigoplus_{j\ge1}C_j/p^rC_j.
\]
To see that its values have finite support, choose $m$ with $p^mx=0$.
For $j\ge m+r$, the coordinate $x_j$ belongs to
$p^{j-m}C_j\subseteq p^rC_j$. The map is surjective because
finitely supported representatives belong to $P_p$.
Clearly $p^rP_p\subseteq\ker\rho_r$. Conversely, suppose that
$x\in\ker\rho_r$ and $p^mx=0$. Choose $y_j\in C_j$ with
$p^ry_j=x_j$ for each $j$. Then $p^{m+r}y_j=0$ for all $j$,
so $y=(y_j)_j\in P_p$ and $x=p^ry$. Thus
$\ker\rho_r=p^rP_p$. Since
$C_j/p^rC_j\cong\Z/p^{\min(j,r)}\Z$, this proves
\eqref{eq:completion-quotients}.

\emph{Step 3: apply the obstruction and realization criterion.}
In particular, $P_p/pP_p$ is countable, while $P_p$ is uncountable
and $p$-adically separated. Theorem~\ref{thm:cardinality} excludes
a half-factorial generating subset. The final assertion follows
from the classical realization criterion,
Theorem~\ref{thm:realization}.
\end{proof}

The argument can be strengthened from nonexistence to the production
of a specific kind of finite obstruction inside every generating set.

\begin{corollary}\label{cor:universal-large-atom}
For every prime $p$ and every generating subset
$S\subseteq P_p\setminus\{0\}$, there is an atom of $\B(S)$
whose cross number is greater than one. More precisely, there are
distinct $s,t\in S$ of a common order $q$ and a nonempty sequence
$C$ over $S\setminus\{s,t\}$ such that $st^{q-1}C$ is an atom.
\end{corollary}

\begin{proof}
A generating set of the uncountable group $P_p$ is uncountable.
Partition it using the sets $S_{n,r}$ in
\eqref{eq:generator-bound}, without assuming half-factoriality.
Since every quotient $P_p/p^rP_p$ is countable, some $S_{n,r}$
contains distinct $s,t$ with $s-t\in p^rP_p$.
Choose a finite $S_0\subseteq S$ containing these elements and
generating a $p^r$-root of $s-t$. Apply
Lemma~\ref{lem:finite-separator} to $\langle S_0\rangle$.
Exactly as in the proof of Theorem~\ref{thm:generator-separation},
the resulting subgroup $K$ gives the hypotheses of
Lemma~\ref{lem:collision}, with $q=p^n$. The correction sequence
there lies in $\F(A)$ with $A\subseteq S_0\cap K$, whereas
$s,t\notin K$. The resulting atom is also an atom over $S$.
\end{proof}

\subsection{Sharpness and smaller counterexamples}

The known countable case gives a converse to the cardinal obstruction
when the quotient modulo $p$ is countable.

\begin{corollary}\label{cor:countable-criterion}
Let $p$ be a prime and let $P$ be a $p$-adically separated abelian
$p$-group with countable $P/pP$. Then $P$ has a half-factorial
generating subset if and only if $P$ is countable.
\end{corollary}

\begin{proof}
Necessity follows from Theorem~\ref{thm:cardinality}. Sufficiency is
the countable-group theorem of Geroldinger--G\"obel
\cite[Theorem~5.3]{GeroldingerGoebel}.
\end{proof}

Both the size of the quotient and the separatedness hypothesis matter.
The following examples make those boundaries explicit.

\begin{example}\label{ex:cardinality-boundaries}
Let $p$ be a prime and $\kappa$ an infinite cardinal. The elementary
abelian group $E=\bigoplus_{\kappa}C_p$ is $p$-adically separated
and has a half-factorial generating set consisting of its cyclic
basis elements. It satisfies $|E|=|E/pE|=\kappa$, so the bound in
Theorem~\ref{thm:cardinality} is attained at every infinite cardinal.

In contrast, for uncountable $\kappa$, the divisible group
$D=\bigoplus_{\kappa}C_{p^\infty}$ belongs to $\HF$ by
Lemma~\ref{lem:directsums} and Proposition~\ref{prop:prufer}, but
$pD=D$ and $|D|=\kappa>\aleph_0$. Thus the separatedness assumption
cannot be omitted. The theorem does not claim that its cardinal
bound is sufficient in general.
\end{example}

The continuum-sized example contains pure counterexamples of every
smaller uncountable cardinality. We include the subgroup construction
rather than assuming that half-factorial realizability is inherited
by subgroups.

\begin{corollary}\label{cor:aleph-one}
Let $p$ be a prime and let $\aleph_1\le\kappa\le2^{\aleph_0}$.
There is a pure subgroup $H_\kappa$ of
$P_p=t(\prod_{j\ge1}\Z/p^j\Z)$ such that
\[
 |H_\kappa|=\kappa,\qquad
 H_\kappa/p^rH_\kappa\cong
 \bigoplus_{j\ge1}\Z/p^{\min(j,r)}\Z\quad(r\ge1).
\]
The group $H_\kappa$ is reduced, unbounded, and $p$-adically
separated, and it has no half-factorial generating subset.
Consequently the least cardinality of a group not realizable as
the ideal class group of a half-factorial Dedekind domain is
$\aleph_1$.
\end{corollary}

\begin{proof}
Put $B=\bigoplus_{j\ge1}\Z/p^j\Z\subseteq P_p$ and choose
$X\subseteq P_p$ of cardinality $\kappa$. Let
$H_0=\langle B\cup X\rangle$. Having constructed $H_i$, choose,
for each $h\in H_i\cap p^rP_p$ and $r\ge1$, one element
$y_{h,r}\in P_p$ with $p^ry_{h,r}=h$, and let $H_{i+1}$ be the
subgroup generated by $H_i$ and all these roots. Each $H_i$ has
cardinality $\kappa$. Hence $H_\kappa=\bigcup_{i\ge0}H_i$ has
cardinality $\kappa$ and satisfies
$H_\kappa\cap p^rP_p=p^rH_\kappa$ for every $r$. Multiplication
by integers prime to $p$ is an automorphism of every $p$-group,
so $H_\kappa$ is pure in $P_p$.

Purity gives an injection
$H_\kappa/p^rH_\kappa\to P_p/p^rP_p$. It is surjective because
$B\subseteq H_\kappa$ and the coordinate calculation in
Theorem~\ref{thm:counterexample} shows that $B$ maps onto
$P_p/p^rP_p$. This proves the displayed quotient formula.
Separatedness is inherited from $P_p$, and it implies reducedness;
unboundedness follows from $B\subseteq H_\kappa$. Theorem~\ref{thm:cardinality}
now excludes half-factorial generators. Finally, every countable
abelian group belongs to $\HF$, whereas $H_{\aleph_1}$ does not.
Apply Theorem~\ref{thm:realization}.
\end{proof}

The nonrealization proof is now complete. It depends only on the
classical criteria in Section~\ref{sec:preliminaries}, finite
abelian-group structure, and the arguments in this section. The
following sections study complementary phenomena and are not
needed for this proof.

\section{Lifting and extensions}\label{sec:extensions}

We compare constructions that preserve half-factoriality with
extensions that cannot do so. Free kernels admit an explicit
block-monoid decomposition, while order-preserving lifts across
torsion kernels impose a splitting condition. We also recall the
classical primary reduction and show that the class $\HF$ is not
closed under pure extensions.

\subsection{An arbitrary free kernel}

The essential observation is that lifts can be moved into a common
nonnegative coordinate cone along a free kernel. The resulting
correction terms are then additive and nonnegative. When the free
subgroup has torsion quotient, the block-monoid decomposition below
is \cite[Proposition~3.7]{GeroldingerGoebel}. The proof given here
allows a quotient with elements of infinite order as well.

\begin{theorem}\label{thm:lifting}
Let $G$ be an abelian group and let
$F=\bigoplus_{i\in I}\Z e_i$ be a free subgroup. Write
$Q=G/F$, and suppose that $T\subseteq Q\setminus\{0\}$ satisfies
$[T]=Q$. There is a subset $S\subseteq G\setminus\{0\}$ with
$[S]=G$ and a monoid isomorphism
\begin{equation}\label{eq:mainiso}
 \B(T)\times\Nzero^{(I)}\ \cong\ \B(S).
\end{equation}
Here $\Nzero^{(I)}$ denotes the monoid of finitely supported
nonnegative integer vectors, with coordinatewise addition.
Consequently, $G/F\in\HF$ implies $G\in\HF$. No purity or
splitting hypothesis on $F\subseteq G$ is required.
\end{theorem}

\begin{proof}
\emph{Step 1: choose nonnegative representatives.}
Since tensoring with $\Q$ is exact,
$\Q\otimes_{\Z}F$ embeds in $\Q\otimes_{\Z}G$. Choose a
$\Q$-linear projection onto this subspace. Composing the natural
map from $G$ with this projection gives a group homomorphism
\[
 \lambda:G\longrightarrow\Q^{(I)},\qquad
 \lambda(e_i)=\varepsilon_i,
\]
where $\varepsilon_i$ is the $i$th coordinate vector. In particular,
$\lambda$ is injective on $F$ and identifies $F$ with
$\Z^{(I)}$.

For each $t\in T$, choose a representative $y_t\in G$. All vectors
$\lambda(y_t)$ have finite support. Adding an appropriate integer
combination of the $e_i$ therefore produces a representative
$x_t\in y_t+F$ such that $\lambda_i(x_t)\ge0$ for every $i$.
Define
\begin{equation}\label{eq:S}
 S=\{x_t:t\in T\}\ \cup\ \{e_i,-e_i:i\in I\}.
\end{equation}
The three parts of this union are disjoint: the $x_t$ lie in
nonzero, pairwise distinct cosets of $F$, and the basis elements
have infinite order. Projecting an arbitrary $g\in G$ to $Q$
expresses its coset as a finite sum of members of $T$. Subtracting
the corresponding $x_t$ leaves an element of $F$, which is a
nonnegative sum of the $e_i$ and $-e_i$. Thus $[S]=G$.

\emph{Step 2: define additive correction terms.}
For $B\in\B(T)$, the element $\sum_t v_t(B)x_t$ lies in $F$.
Write
\begin{equation}\label{eq:correction}
 \sum_{t\in T}v_t(B)x_t=\sum_{i\in I}a_i(B)e_i.
\end{equation}
The coefficients $a_i(B)$ are integers with finite support.
Applying $\lambda$ shows that
$a_i(B)=\sum_t v_t(B)\lambda_i(x_t)\ge0$. The map
$B\mapsto(a_i(B))_i$ is additive under concatenation. Hence
\begin{equation}\label{eq:hat}
 \widehat B=
 \prod_{t\in T}x_t^{v_t(B)}
 \prod_{i\in I}(-e_i)^{a_i(B)}
 \quad\text{belongs to }\B(S).
\end{equation}
Define
\begin{equation}\label{eq:Phi}
 \Phi\bigl(B,(n_i)_i\bigr)=
 \widehat B\prod_{i\in I}\bigl(e_i(-e_i)\bigr)^{n_i}.
\end{equation}
Additivity of the correction vector shows that $\Phi$ is a monoid
homomorphism.

\emph{Step 3: recover a block from its unique components.}
To prove bijectivity, let $C\in\B(S)$. Its sequence expression is
uniquely of the form
\[
 C=\prod_{t\in T}x_t^{v_t}
   \prod_{i\in I}e_i^{u_i}(-e_i)^{w_i}.
\]
Projection to $Q$ shows that $B=\prod_t t^{v_t}$ belongs to
$\B(T)$. The zero-sum condition in $G$, together with
\eqref{eq:correction} and independence of the $e_i$, gives
$w_i=a_i(B)+u_i$ for every $i$. Therefore
$C=\Phi(B,(u_i)_i)$. The lift multiplicities determine $B$, and
the $e_i$-multiplicities determine $(u_i)_i$, so this expression is
unique. This proves \eqref{eq:mainiso}. A product with a free
commutative monoid is half-factorial whenever the other factor is
half-factorial, proving the last assertion.
\end{proof}

The same isomorphism identifies the entire factorization-length
structure, not only half-factoriality.

\begin{corollary}\label{cor:lengths}
Let $G$ be an abelian group, let $F=\bigoplus_{i\in I}\Z e_i$
be a free subgroup, and let $T\subseteq(G/F)\setminus\{0\}$
satisfy $[T]=G/F$. Choose the representatives $x_t$, the set $S$,
and the isomorphism $\Phi$ by the construction in the proof of
Theorem~\ref{thm:lifting}. Then the atoms of $\B(S)$ are precisely the sequences $\widehat U$, with
$U\in\A(\B(T))$, and the pairs $e_i(-e_i)$, with $i\in I$.
For every $B\in\B(T)$ and $(n_i)_i\in\Nzero^{(I)}$,
\begin{equation}\label{eq:lengths}
 \LL_{\B(S)}\!\left(\Phi\bigl(B,(n_i)_i\bigr)\right)
 =\LL_{\B(T)}(B)+\sum_{i\in I}n_i.
\end{equation}
In particular, if $T$ is half-factorial and its factorization-length
function is $\ell_T$, then the factorization-length function on $\B(S)$
is given by
$\ell_S(\Phi(B,(n_i)_i))=\ell_T(B)+\sum_i n_i$.
\end{corollary}

\begin{proof}
An atom of a direct product of reduced atomic monoids has exactly
one nonidentity component, and that component is an atom. In the
factor $\Nzero^{(I)}$, the atoms are the coordinate vectors and
the factorization length of $(n_i)_i$ is $\sum_i n_i$. Transport
these statements through $\Phi$.
\end{proof}

The construction applies to nonsplit extensions in a concrete way.

\begin{example}\label{ex:localization}
Let $p$ be a prime, $G=\Z[1/p]$, and $F=\Z$. Then
$G/F\cong C_{p^\infty}$. For $c_n=p^{-n}+\Z$, take
$T=\{c_n:n\ge1\}$ and the representatives $x_{c_n}=p^{-n}$.
Theorem~\ref{thm:lifting} gives
\[
 S=\{-1,1\}\cup\{p^{-n}:n\ge1\},\qquad
 \B(S)\cong\B(T)\times\Nzero.
\]
The half-factoriality of $T$ is verified in
Proposition~\ref{prop:prufer}. A block $B\in\B(T)$ has integral
weight $m=\sum_n v_{c_n}(B)p^{-n}$, and its lift is
$\widehat B=\prod_n(p^{-n})^{v_{c_n}(B)}(-1)^m$.
The exact sequence does not split: a nonzero torsion quotient
$C_{p^\infty}$ cannot be a subgroup of the torsion-free group
$\Z[1/p]$. Also, $F$ is not pure in $G$, since $1\in pG\cap F$
but $1\notin pF$.
\end{example}

\subsection{The classical primary reduction and positive consequences}

The reduction from arbitrary groups to reduced uncountable unbounded
primary groups is due to Geroldinger--G\"obel
\cite[pp.~855--856]{GeroldingerGoebel}. We give the argument because
it explains the role of free kernels and distinguishes the known
reduction from the separation theorem of Section~\ref{sec:cardinality}.

\begin{theorem}[Geroldinger--G\"obel reduction]\label{thm:reduction}
The following assertions are equivalent.
\begin{enumerate}[label=\textup{(\arabic*)},leftmargin=2.2em]
\item Every abelian group belongs to $\HF$.
\item Every torsion abelian group belongs to $\HF$.
\item For every prime $p$, every abelian $p$-group belongs to $\HF$.
\item For every prime $p$, every reduced abelian $p$-group belongs
      to $\HF$.
\item For every prime $p$, every reduced uncountable unbounded
      abelian $p$-group belongs to $\HF$.
\end{enumerate}
Consequently each assertion is equivalent to an affirmative answer
to the universal half-factorial Dedekind realization problem.
\end{theorem}

\begin{proof}
The implications
$(1)\Rightarrow(2)\Rightarrow(3)\Rightarrow(4)\Rightarrow(5)$
are immediate. Assertion (5) implies (4) by the countable and bounded cases recalled in
Section~\ref{sec:preliminaries}. We prove the other converse steps.

Assume (4), and let $P$ be a $p$-group. Its maximal divisible
subgroup $D$ is a direct summand, so $P=D\oplus R$ with $R$
reduced. Moreover, $D$ is a direct sum of copies of
$C_{p^\infty}$. These are standard structure results for abelian
groups \cite{Fuchs}. By Proposition~\ref{prop:prufer},
Lemma~\ref{lem:directsums}, and (4), both summands and hence $P$
belong to $\HF$. This proves (3).

Assume (3). Every torsion abelian group is the direct sum of its
$p$-primary components \cite{Fuchs}. Apply
Lemma~\ref{lem:directsums} to obtain (2).

Finally, assume (2), and let $G$ be any abelian group. Choose a
maximal $\Z$-linearly independent subset $\{e_i:i\in I\}$ of
$G$ and let $F=\bigoplus_i\Z e_i$. For every $g\in G$, either
$g\in F$ or maximality gives a relation
$ng+\sum_i n_i e_i=0$ with $n\ne0$: a relation with $n=0$
would contradict independence of the $e_i$. Thus a nonzero
multiple of every coset in $G/F$ is zero, and $G/F$ is torsion.
By (2), $G/F\in\HF$. Theorem~\ref{thm:lifting} gives $G\in\HF$,
proving (1). The final assertion follows from
Theorem~\ref{thm:realization}.
\end{proof}

The reduction is a logical equivalence between universal assertions;
Theorem~\ref{thm:counterexample} shows that all of them are false.
For an individual quotient, however, the lifting theorem remains
a useful positive construction.

\begin{corollary}\label{cor:classes}
Let $G$ be an abelian group and let $F$ be a free subgroup. If
$G/F$ is countable, divisible, or a direct sum of cyclic groups
and Pr\"ufer groups, then $G\in\HF$. In particular, every
torsion-free abelian group of finite rank is the ideal class group
of a half-factorial Dedekind domain.
\end{corollary}

\begin{proof}
The countable case is \cite[Theorem~5.3]{GeroldingerGoebel}.
A divisible group is a direct sum of copies of $\Q$ and of
Pr\"ufer groups \cite{Fuchs}; these summands belong to $\HF$ by
the countable case. The cyclic and Pr\"ufer case follows from
Lemma~\ref{lem:directsums} and Proposition~\ref{prop:prufer}.
Apply Theorem~\ref{thm:lifting}. A torsion-free group of finite
rank embeds in $\Q^r$ for some finite $r$, and is therefore
countable. Its realization follows from the countable case and
Theorem~\ref{thm:realization}.
\end{proof}

More generally, the quotient may be any group already known to
belong to $\HF$. In particular, the Warfield-group theorem
\cite[Theorem~5.1]{GeroldingerGoebel} applies; the corresponding
full-rank quotient consequence is explicitly
\cite[Corollary~5.2]{GeroldingerGoebel}. 

\subsection{Root adjunction and torsion kernels}

A single root adjunction admits an exact criterion. The special
case in which the lower multiple is already among the generators
is \cite[Proposition~3.3]{GeroldingerGoebel}. The formulation below
also covers a lower multiple not originally in the generating set.

\begin{proposition}\label{prop:root}
Let $G$ be a torsion abelian group, let $H\subseteq G$ be a subgroup,
and let $T\subseteq H\setminus\{0\}$ generate $H$. For
$x\in G\setminus H$, put $m=\ord(x+H)$ and $h=mx\in H$.
There is a surjective homomorphism
\[
 \theta:\B(T\cup\{x\})\longrightarrow\B(T\cup\{h\}),
 \qquad A x^{mk}\longmapsto A h^k,
\]
which preserves all sets of factorization lengths. Consequently,
$T\cup\{x\}$ is half-factorial if and only if $T\cup\{h\}$ is
half-factorial. The statement includes $h=0$.
\end{proposition}

\begin{proof}
The multiplicity of $x$ in a block is divisible by $m$, by
projection to $G/H$, and substitution of $x^m$ by $h$ preserves
its sum. If $h\notin T$, the displayed map is an isomorphism:
replace every occurrence of $h$ by $x^m$ to obtain its inverse.

If $h\in T$, distinguish temporarily the occurrences of $h$
coming from $x^m$ from the original occurrences of $h$ in $A$.
A factorization of $\theta(Ax^{mk})$ can be lifted by distributing
these marked occurrences among its factors, replacing each marked
occurrence by $x^m$. The resulting factors are nonempty blocks
and multiply to $Ax^{mk}$. Conversely, the image of a nonempty
block is nonempty; if the image of an atom factored, this lifting
would factor the atom. Surjectivity follows by taking $k=0$.
Therefore the map preserves exactly the sets of lengths.
\end{proof}

For a sequence over a subgroup or a quotient we compute its cross
number using the orders in that group, indicating the group by a
subscript when necessary.

\begin{lemma}\label{lem:defect}
Let
\[
 0\longrightarrow K\longrightarrow G
   \mathrel{\mathop{\longrightarrow}^{\pi}} Q
   \longrightarrow0
\]
be an exact sequence of torsion abelian groups. Let
$A\subseteq K\setminus\{0\}$ and
$T\subseteq Q\setminus\{0\}$ be half-factorial generating sets.
Choose $x_t\in G$ with $\pi(x_t)=t$ and put
$S=A\cup\{x_t:t\in T\}$. For
$U\in\A(\B(T))$, define
\[
 V_U=\prod_{t\in T}x_t^{v_t(U)},\qquad
 d(U)=1-\sum_{t\in T}\frac{v_t(U)}{\ord_G(x_t)}.
\]
If $S$ is half-factorial, then there is a sequence $C\in\F(A)$
such that $V_UC$ is an atom of $\B(S)$ and
\begin{equation}\label{eq:defect}
 \sigma(C)=-\sigma(V_U),\qquad \mathsf k_K(C)=d(U).
\end{equation}
If $\sigma(V_U)\ne0$, then $d(U)>0$ and at least one letter
$t$ in $U$ satisfies $\ord_G(x_t)>\ord_Q(t)$.
\end{lemma}

\begin{proof}
Since $\sigma(V_U)\in K$ and $A$ generates $K$ as a monoid,
choose $D\in\F(A)$ with $\sigma(D)=-\sigma(V_U)$. Factor the
block $V_UD$ into atoms over $S$. Projection to $Q$, omitting
the letters in $A$, sends each factor to a block dividing $U$.
Because $U$ is an atom, exactly one factor contains all the
letters of $V_U$; all other factors contain only letters of
$A$. Write that distinguished factor as $V_UC$, with $C\mid D$.
Half-factoriality and Proposition~\ref{prop:cross} give
$\mathsf k_G(V_U)+\mathsf k_K(C)=1$.

If $\sigma(V_U)\ne0$, then $C$ is nonempty and its cross number
is positive. Also $\ord_Q(t)$ divides $\ord_G(x_t)$, while
$\mathsf k_Q(U)=1$. Thus $d(U)>0$ requires strict increase of
the order of at least one lifted letter.
\end{proof}

In particular, there is no room for a nonempty correction when all
lifted orders equal the corresponding quotient orders. This yields
an equivalence, rather than merely a necessary numerical test.

\begin{theorem}\label{thm:splitting}
Let $0\to K\to G\xrightarrow{\pi}Q\to0$ be an exact sequence
of torsion abelian groups. Let $A\subseteq K\setminus\{0\}$
and $T\subseteq Q\setminus\{0\}$ be half-factorial generating
sets. Choose representatives $x_t\in G$ such that
\[
 \pi(x_t)=t,\qquad \ord_G(x_t)=\ord_Q(t)\quad(t\in T),
\]
and put $S=A\cup\{x_t:t\in T\}$. The following are equivalent.
\begin{enumerate}[label=\textup{(\arabic*)},leftmargin=2.2em]
\item $S$ is half-factorial.
\item For every $U\in\A(\B(T))$, one has
      $\sum_t v_t(U)x_t=0$ in $G$.
\item The assignment $t\mapsto x_t$ extends to a homomorphism
      $s:Q\to G$ with $\pi s=\operatorname{id}_Q$.
\end{enumerate}
When these conditions hold, $G=K\oplus s(Q)$ and
$\B(S)\cong\B(A)\times\B(T)$. In particular, if the exact
sequence does not split, no choice of order-preserving
representatives yields a half-factorial set of the displayed form.
\end{theorem}

\begin{proof}
Under the order hypothesis, $d(U)=1-\mathsf k_Q(U)=0$ for every
quotient atom $U$. Lemma~\ref{lem:defect} proves $(1)\Rightarrow(2)$.

Assume (2). Since every block over $T$ factors into atoms, the
lifted sum of every such block is zero. Consider a finite integer
relation $\sum_t n_t t=0$ in $Q$. For each occurring $t$, choose
$m_t\ge0$ with $m_t\equiv n_t\pmod{\ord_Q(t)}$. The sequence
$\prod_t t^{m_t}$ is a block, so $\sum_t m_tx_t=0$ by the
preceding observation. Since
$\ord_G(x_t)=\ord_Q(t)$, replacing $m_t$ by $n_t$ does not
change its lifted sum. Thus $\sum_t n_tx_t=0$.

All defining group relations among $T$ are therefore respected.
The assignment extends to a homomorphism $s:Q\to G$, and
$\pi s=\operatorname{id}_Q$ because this equality holds on the
generating set $T$. This proves $(2)\Rightarrow(3)$.

Finally, a section gives $G=K\oplus s(Q)$. The two parts of $S$
lie in the two summands, so the block monoid is the displayed
product. Both factors are half-factorial. This proves
$(3)\Rightarrow(1)$ and the remaining assertions.
\end{proof}

The hypotheses concern the chosen generating sets. They do not say
that a nonsplit torsion extension fails to belong to $\HF$.
The following countable example makes that limitation explicit.

\begin{example}\label{ex:positive-nonsplit}
Let $p$ be a prime and let $E=\bigoplus_{n\ge1}\langle e_n\rangle$ with
$\ord(e_n)=p^n$, and define an epimorphism
$\rho:E\to C_{p^\infty}$ by
$\rho(e_n)=c_n=1/p^n+\Z$. Put $k_n=e_n-pe_{n+1}$.
Then
\[
 K=\ker\rho=\bigoplus_{n\ge1}\langle k_n\rangle,
 \qquad \ord(k_n)=p^n.
\]
Indeed, a relation among finitely many $k_n$ is trivial by
successively inspecting its highest $e_{n+1}$-coordinate. To
identify the kernel, replace $e_n$ successively by
$k_n+pe_{n+1}$ in a finitely supported element of $\ker\rho$;
the remaining last coefficient is divisible by the order of its
basis element and hence contributes zero.

Both $A=\{k_n:n\ge1\}$ and $T=\{c_n:n\ge1\}$ are
half-factorial generating sets. The representatives $x_{c_n}=e_n$
preserve orders, but the exact sequence does not split, since
$E$ is reduced. Thus $A\cup\{e_n:n\ge1\}$ is not
half-factorial. Directly, for every $n\ge1$ the sequence
$e_{n+1}^{p}k_ne_n^{p^n-1}$ is an atom with cross number
$1+1/p^n$: a zero-sum subsequence must either omit $k_n$ and be
empty, or contain $k_n$ and then have all the indicated
multiplicities. Nevertheless, $\{e_n:n\ge1\}$ alone is a
half-factorial generating set of $E$. Thus the obstruction comes
from insisting on adjoining the prescribed kernel generators.
\end{example}

\subsection{The pure extension associated to the torsion-completion}

The group from Theorem~\ref{thm:counterexample} also exhibits a
failure of extension closure. We identify its pure direct-sum
subgroup and divisible quotient explicitly. By the
\emph{torsion-completion} of a $p$-group $B$ we mean the torsion
subgroup of its $p$-adic completion
$\widehat B=\varprojlim_r B/p^rB$; see \cite{FuchsII}.

\begin{proposition}\label{prop:completion}
Fix a prime $p$. Set $C_n=\Z/p^n\Z$,
$B=\bigoplus_{n\ge1}C_n$, and
\begin{equation}\label{eq:completion}
 P=\left\{(a_n)_n\in\prod_{n\ge1}C_n:
       p^r a_n=0\text{ for every }n\text{ and some }r\ge0\right\}.
\end{equation}
Then $P$ is the torsion-completion of $B$, is reduced, uncountable,
and unbounded, and $B$ is pure in $P$. The quotient $D=P/B$ is a
nonzero divisible $p$-group.

Let $A\subseteq B\setminus\{0\}$ and
$T\subseteq D\setminus\{0\}$ be any half-factorial generating
sets. Order-preserving representatives $x_t\in P$ of all $t\in T$
exist, but for every such choice the set
$A\cup\{x_t:t\in T\}$ is not half-factorial.
\end{proposition}

\begin{proof}
The $p$-adic completion of $B$ identifies with the subgroup of
$\prod_n C_n$ consisting of sequences whose coordinates belong to
$p^kC_n$ for all but finitely many $n$, for every $k$.
This follows by viewing $B/p^kB$ as the direct sum of the
coordinate quotients $C_n/p^kC_n$ and taking the inverse limit.
If a sequence is annihilated by $p^r$, its $n$th coordinate for
$n>r$ belongs to $p^{n-r}C_n$. Thus every sequence in
\eqref{eq:completion} lies in the completion. Multiplication by an integer prime to $p$ is invertible in every
coordinate. Thus every torsion element of the completion is killed
by a power of $p$, and its torsion subgroup is exactly $P$.

For each fixed coordinate $n$, every element of $p^kP$ has zero
$n$th coordinate when $k\ge n$. Hence
$\bigcap_{k\ge0}p^kP=0$, which implies that $P$ is reduced.
Moreover, $P[p]=\prod_{n\ge1}C_n[p]$ is a product of countably
many nontrivial groups of order $p$, so $P$ is uncountable.
It is unbounded because it contains $B$.

If $b\in B\cap p^kP$, choose $y\in P$ with $p^ky=b$ and
truncate $y$ to the finite support of $b$. The truncation lies in
$B$ and still maps to $b$, proving $B\cap p^kP=p^kB$.
Multiplication by any integer prime to $p$ is an automorphism on
both groups, so $B$ is pure in $P$.

To prove divisibility of $D$, let $x=(x_n)_n\in P$ be annihilated
by $p^r$, with $r\ge1$. Truncate the first $r$ coordinates to
obtain $b\in B$. For $n>r$, the coordinate $x_n$ is divisible
by $p$ in $C_n$. Choose roots coordinatewise to obtain $y$ with
$py=x-b$. The roots are annihilated by $p^{r+1}$, so $y\in P$.
It follows that multiplication by $p$ on $D$ is surjective, and
$D$ is divisible since it is a $p$-group. The element
$(p^{n-1}+p^n\Z)_{n\ge1}$ belongs to $P\setminus B$, so
$D\ne0$.

The group $B$ has a half-factorial generating set by its cyclic
decomposition; $D$ has one by its decomposition into Pr\"ufer
groups and Proposition~\ref{prop:prufer}. For any element
$t\in D$ of order $p^r$, choose a representative $y\in P$.
Purity gives
$p^ry\in B\cap p^rP=p^rB$. Choose $b\in B$ with
$p^rb=p^ry$. Then $x_t=y-b$ represents $t$ and has order exactly
$p^r$.

Finally, $0\to B\to P\to D\to0$ cannot split: a section would
embed the nonzero divisible group $D$ in the reduced group $P$.
The last assertion therefore follows from
Theorem~\ref{thm:splitting}.
\end{proof}

\begin{remark}\label{rem:completion-limit}
Proposition~\ref{prop:completion}, considered on its own, excludes
only a specified class of lifting constructions. A putative
half-factorial generating set of $P$ need not contain generators
of this particular subgroup $B$, have a half-factorial image
in $D$, or preserve orders under projection. These possibilities
are all excluded separately by Theorem~\ref{thm:counterexample},
whose proof makes none of those assumptions. Example~\ref{ex:positive-nonsplit}
shows that the first limitation matters even in a countable group
already known to belong to $\HF$. Theorem~\ref{thm:weak-existence}
also applies to $P$, so a weakly half-factorial generating set does
exist in this test case.
\end{remark}

The pure extension just computed yields a closure failure for
existence of half-factorial generators, not merely for a selected
choice of representatives.

\begin{corollary}\label{cor:pure-extension-failure}
The class $\HF$ is not closed under extensions, even under pure
extensions of abelian $p$-groups. For every prime $p$, the exact
sequence
\[
 0\longrightarrow B\longrightarrow P_p\longrightarrow P_p/B
 \longrightarrow0,\qquad B=\bigoplus_{j\ge1}C_{p^j},
\]
has pure kernel, $B\in\HF$, and $P_p/B\in\HF$, but
$P_p\notin\HF$.
\end{corollary}

\begin{proof}
Proposition~\ref{prop:completion} proves purity and divisibility of
$P_p/B$. Direct sums of cyclic groups and divisible groups belong
to $\HF$, whereas Theorem~\ref{thm:counterexample} excludes $P_p$.
\end{proof}

\section{Characters and compatibility of generating sets}\label{sec:characters}

Integrality of cross numbers is weaker than half-factoriality.
We first show that every torsion abelian group admits generators
satisfying this weaker condition. We then impose a prescribed
character on each bounded primary group and examine the obstructions
to extending or thinning these generators. The completion examples
show that these positive constructions cannot be made globally
half-factorial in general.

\subsection{Weak half-factoriality and characters}

A subset $T$ of a torsion group is \emph{weakly half-factorial} if
every atom of $\B(T)$ has integral cross number. Since every block
factors into atoms, this is equivalent to integrality of the cross
number of every block. The terminology and finite-group character
criterion are standard
\cite[Lemma~4.1]{GaoGeroldinger}; see also
\cite{RadziejewskiSchmid,Schmid2005}. The next proof works without
a boundedness or finiteness assumption.

\begin{lemma}\label{lem:character}
Let $G$ be a torsion abelian group and let
$T\subseteq G\setminus\{0\}$ generate $G$. The following are
equivalent: every $B\in\B(T)$ has integral cross number; and there
is a homomorphism $\chi:G\to\Q/\Z$ satisfying
$\chi(t)=1/\ord(t)+\Z$ for every $t\in T$.
\end{lemma}

\begin{proof}
A character with the stated values gives
$\mathsf k(B)+\Z=\chi(\sigma(B))=0$ for each block $B$.
Conversely, assign $t$ the value $1/\ord(t)+\Z$. Given a finite
integer relation $\sum_t n_t t=0$, replace each $n_t$ by a
nonnegative integer $m_t$ congruent to $n_t$ modulo $\ord(t)$.
Then $\prod_t t^{m_t}$ is a block, so
$\sum_t m_t/\ord(t)$ is integral. The same is true of
$\sum_t n_t/\ord(t)$. Thus the assignment respects every relation
and defines the required homomorphism.
\end{proof}

The character condition is always achievable when the generating
set is allowed to vary. The following proof uses only a linear
functional on the socle and the injectivity of divisible abelian
groups. The finite-group character theory is developed in
\cite{RadziejewskiSchmid}; the following construction treats
arbitrary cardinality and unbounded exponent.

\begin{theorem}\label{thm:weak-existence}
Every torsion abelian group $G$ has a generating subset
$T\subseteq G\setminus\{0\}$ and a homomorphism
$\chi:G\to\Q/\Z$ such that
\begin{equation}\label{eq:weak-character}
 \chi(t)=\frac{1}{\ord(t)}+\Z\qquad(t\in T).
\end{equation}
Consequently every zero-sum sequence over $T$ has integral cross
number. In particular, every reduced abelian $p$-group has a weakly
half-factorial generating subset.
\end{theorem}

\begin{proof}
We first treat a nonzero abelian $p$-group $P$. Write
$P[p]=\{x\in P:px=0\}$ and identify the subgroup of order $p$ in
$C_{p^\infty}\subseteq\Q/\Z$ with the additive group of
$\mathbb F_p$, sending $1$ to $1/p+\Z$.

\emph{Step 1: bounded primary groups.}
Suppose first that $P$ has exponent $p^N$. Choose $z\in P$ of
order $p^N$ and extend the map
$z\mapsto 1/p^N+\Z$ on $\langle z\rangle$ to a homomorphism
$\chi:P\to C_{p^\infty}$. This extension exists because divisible
abelian groups are injective \cite{Fuchs}. Put
\[
 T=\{t\in P\setminus\{0\}:\chi(t)=1/\ord(t)+\Z\}.
\]
Then $z\in T$. For $g\in P$, write
$\chi(g)=a/p^N+\Z$ with $0\le a<p^N$. If $p\nmid a$, then
$g$ has order $p^N$ and multiplication by an inverse of $a$
modulo $p^N$ produces a member of $T$ generating the same cyclic
subgroup as $g$. If $p\mid a$, apply the same argument to $g+z$,
whose character value has numerator $a+1$ prime to $p$.
Then $g=(g+z)-z\in\langle T\rangle$. Thus $T$ generates $P$.

\emph{Step 2: a functional on the socle in the unbounded case.}
Now suppose that $P$ is unbounded. For $n\ge0$ let
$V_n=(p^nP)[p]\subseteq P[p]$. These are nonzero descending
$\mathbb F_p$-subspaces. There is a linear functional
$f:P[p]\to\mathbb F_p$ that is nonzero on every $V_n$.
Indeed, if $\bigcap_n V_n\ne0$, prescribe $f(v)=1$ on a nonzero
vector in the intersection and extend to a basis. If the intersection
is zero, choose increasing integers $n_j$ with
$V_{n_j}\supsetneq V_{n_{j+1}}$ and vectors
$v_j\in V_{n_j}\setminus V_{n_{j+1}}$. The $n_j$ are unbounded,
and the $v_j$ are linearly independent: in a nontrivial finite
relation, the term with the least index cannot lie in the subspace
containing all later terms. Prescribe $f(v_j)=1$ and extend to a
basis. This again gives the required functional.

\emph{Step 3: extend the functional and normalize generators.}
Regard $f$ as a homomorphism to the subgroup of order $p$ in
$C_{p^\infty}$ and extend it to
$\chi:P\to C_{p^\infty}$. For every $n\ge1$, choose
$v\in V_{n-1}$ with $f(v)\ne0$ and write $v=p^{n-1}y$.
Then $y$ and $\chi(y)$ both have order $p^n$. Multiplying $y$
by a suitable integer prime to $p$ gives an element $z_n$ of
order $p^n$ with $\chi(z_n)=1/p^n+\Z$.

Define $T$ by the same displayed formula as in the bounded case.
If $g\in P$ is annihilated by $p^m$, then
$\chi(g+z_{m+1})=a/p^{m+1}+\Z$ for an integer
$a\equiv1\pmod p$. Hence $g+z_{m+1}$ has order $p^{m+1}$.
For an inverse $u$ of $a$ modulo $p^{m+1}$, the element
$t=u(g+z_{m+1})$ belongs to $T$, and
$g=at-z_{m+1}\in\langle T\rangle$. Thus $T$ generates $P$.
Since $P$ is torsion, group generation is monoid generation.

\emph{Step 4: combine the primary components.}
For a general torsion group use its primary decomposition
$G=\bigoplus_p G_p$. Apply the construction to each nonzero
$G_p$, take the union of the generating subsets, and define
$\chi((g_p)_p)=\sum_p\chi_p(g_p)$. This sum is finite for every
$g\in G$. Equation~\eqref{eq:weak-character} follows, and
Lemma~\ref{lem:character} gives integrality. For the trivial group
use $T=\varnothing$ and the zero homomorphism.
\end{proof}

For a homomorphism $\chi:P\to\Q/\Z$ on an abelian $p$-group, put
\begin{equation}\label{eq:admissible-set}
 W_\chi(P)=\{x\in P\setminus\{0\}:
                    \chi(x)=1/\ord(x)+\Z\}.
\end{equation}
Every block over this set has integral cross number, whether or not
it generates $P$. We refer to $W_\chi(P)$ as the full character set; the next results
contrast its arithmetic with that of suitable subsets.

\subsection{A prescribed character on every bounded layer}

The bounded-group existence assertion is classical. The next
construction additionally prescribes the character and gives a
uniform sequence length for all block atoms. We use the standard
cyclic decomposition of a bounded abelian $p$-group \cite{Fuchs}.

\begin{theorem}
\label{thm:bounded-character}
Let $p$ be a prime, let $N\ge1$, and let $P$ be an abelian group of
exponent $p^N$. Suppose that $\chi:P\to\Q/\Z$ has image the subgroup
of order $p^N$. Then there is a half-factorial generating subset
$T\subseteq W_\chi(P)$ such that every $t\in T$ has order $p^N$.
Every atom of $\B(T)$ has sequence length $p^N$, and
\[
 \LL_{\B(T)}(B)=\{|B|/p^N\}\qquad(B\in\B(T)).
\]
More explicitly, choose $z\in P$ with
$\chi(z)=1/p^N+\Z$ and write
$\ker\chi=\bigoplus_{i\in I}\langle e_i\rangle$, where
$\ord(e_i)=p^{r_i}$ and $1\le r_i\le N$. One may take
\begin{equation}\label{eq:bounded-shear}
 T=\{z\}\cup\{z+e_i:i\in I\}.
\end{equation}
\end{theorem}

\begin{proof}
The element $z$ has order $p^N$, and $\chi$ restricts to an
isomorphism from $\langle z\rangle$ onto its image. Hence
$P=\langle z\rangle\oplus\ker\chi$. The displayed set generates
$P$, since $e_i=(z+e_i)-z$. Group generation is monoid generation
in a torsion group. Its elements are distinct, have character value
$1/p^N+\Z$, and therefore have order $p^N$.

A sequence $B=z^{a_0}\prod_i(z+e_i)^{a_i}$ is a block precisely
when $p^{r_i}\mid a_i$ for every $i$ and
$p^N\mid a_0+\sum_i a_i$. Partition the occurrences of $z+e_i$
into packets of $p^{r_i}$ identical letters, and regard each
occurrence of $z$ as a packet of size one. The packet sizes are
powers of $p$ not exceeding $p^N$. Whenever their total size is
$mp^N$, they can be grouped into $m$ collections of size $p^N$.
Indeed, the number of packets of the smallest size $p^j<p^N$
is divisible by $p$, by reducing the total size modulo $p^{j+1}$.
Bundle them in groups of $p$ and repeat, retaining the constituent
packets. This ends with $m$ packets of size $p^N$.

Each resulting collection is a zero-sum subsequence: its
$\ker\chi$-coordinates vanish packet by packet, and its
$z$-coordinate has total coefficient $p^N$. Conversely, every
nonempty block has sequence length a positive multiple of $p^N$,
so a block of length $p^N$ is an atom. Thus all atoms have that
length, proving the formula and half-factoriality.
\end{proof}

The same character can be used on all bounded layers of an
unbounded group. The assertion concerns the individual layers;
it does not assert that the resulting sets are nested.

\begin{corollary}\label{cor:bounded-layers}
Let $p$ be a prime and let $P$ be an unbounded abelian $p$-group. There is a homomorphism
$\chi:P\to C_{p^\infty}\subseteq\Q/\Z$ such that, for every
$n\ge1$, the subgroup $P[p^n]=\{x\in P:p^nx=0\}$ has a
half-factorial generating set $T_n\subseteq W_\chi(P)$ consisting
entirely of elements of order $p^n$. Consequently,
$\bigcup_{n\ge1}T_n$ is a weakly half-factorial generating set
of $P$.
\end{corollary}

\begin{proof}
The unbounded-case construction in
Theorem~\ref{thm:weak-existence} gives, for every $n$, an element
$z_n\in P[p^n]$ with $\chi(z_n)=1/p^n+\Z$. Hence the restriction
of $\chi$ to $P[p^n]$ has image of order $p^n$. Apply
Theorem~\ref{thm:bounded-character} to each restriction. The union
generates $P=\bigcup_nP[p^n]$ and is contained in $W_\chi(P)$,
so its blocks have integral cross number.
\end{proof}

The completion example converts this lack of an asserted compatibility
into an actual obstruction, for every prime.

\begin{corollary}\label{cor:increasing-unions}
For every prime $p$, the group $P_p=t(\prod_{j\ge1}C_{p^j})$ is a
countable increasing union of subgroups belonging to $\HF$, but
$P_p\notin\HF$. Moreover, there do not exist nested half-factorial
sets $T_n\subseteq P_p[p^n]\setminus\{0\}$ generating $P_p[p^n]$
for all $n\ge1$. Thus $\HF$ is not closed under countable increasing
unions of subgroups.
\end{corollary}

\begin{proof}
Each $P_p[p^n]$ is bounded, hence belongs to $\HF$, and their union
is $P_p$. Nonmembership follows from
Theorem~\ref{thm:counterexample}. If the sets $T_n$ were nested,
their union would be half-factorial by finite character and would
generate $P_p$, a contradiction.
\end{proof}

Compatibility does not follow merely from using the same character,
even for two successive bounded layers.

\begin{example}\label{ex:incompatible-layers}
Let $P=\langle e\rangle\oplus\langle a\rangle
\oplus\langle b\rangle$, with $\ord(e)=4$ and
$\ord(a)=\ord(b)=2$. Define
$\chi(e)=1/4+\Z$ and $\chi(a)=\chi(b)=1/2+\Z$. Then
\[
 T_1=\{2e,a,b\},\qquad T_2=\{e,a-e,b-e\}
\]
are half-factorial generating sets of $P[2]$ and $P[4]=P$,
respectively, and both lie in $W_\chi(P)$.
The first assertion follows from independence of $2e,a,b$.
For the second, observe that
$\ker\chi=\langle a-2e\rangle\oplus\langle b-2e\rangle$ and
apply Theorem~\ref{thm:bounded-character} with $z=e$.

Nevertheless, their union is not half-factorial. The sequence
\begin{equation}\label{eq:layer-atom}
 U=(2e)\,a\,b\,(a-e)\,(b-e)
\end{equation}
is an atom with cross number two. To check minimality, write the
multiplicities of $2e,a,b,a-e,b-e$ in a subsequence as
$c,\alpha,\beta,u,v\in\{0,1\}$. Its sum vanishes exactly when
$\alpha=u$, $\beta=v$, and $2c-u-v\equiv0\pmod4$.
These conditions force all five multiplicities to be zero or all
to be one. This proves the claim. Thus individually valid
character-compatible constructions do not establish the
half-factoriality of their union.
\end{example}

A different proposed compatibility condition would be to require
$pz_{n+1}=z_n$ for pivots of order $p^n$. In a reduced group such
an infinite sequence cannot exist: its generated subgroup would
be a nonzero copy of $C_{p^\infty}$. This observation excludes
that particular choice of compatible pivots, not all possible
half-factorial generating sets.

\subsection{Finite obstructions to extension and thinning}

The following examples distinguish failure of a chosen generating
set from nonexistence of all half-factorial generating sets. In
particular, a character controls the integrality of cross numbers,
not their value on atoms.

\begin{example}\label{ex:weak}
Let $G=C_2^3$ have basis $a,b,c$, put $d=a+b+c$, and set
$T=\{a,b,c,d\}$. The character that sends each basis element to
$1/2+\Z$ sends every member of $T$ to $1/2+\Z$.
Lemma~\ref{lem:character} therefore applies. Nevertheless,
$U=abcd$ is an atom with $\mathsf k(U)=2$. Indeed, a zero-sum
subsequence not containing $d$ must be empty, while a zero-sum
subsequence containing $d$ must contain all of $a,b,c$.
The two factorizations
\[
 U\cdot U=(a^2)(b^2)(c^2)(d^2)
\]
have lengths two and four. Thus $T$ is not half-factorial.
This is not a counterexample to the universal realization problem: the
basis subset $\{a,b,c\}$ is a half-factorial generating set of
$C_2^3$.
\end{example}

Half-factoriality of a chosen set cannot be transported through
an arbitrary quotient or multiplication map. The next example is
completely explicit.

\begin{proposition}\label{prop:doubling}
There is a finite abelian $2$-group $G$ with a half-factorial
generating subset $S$ such that $2S$ is a non-half-factorial
generating subset of $2G$. Consequently, neither homomorphic
images of half-factorial subsets nor adjoining all their
$2$-multiples necessarily preserves half-factoriality.
\end{proposition}

\begin{proof}
Take $G=C_4^3$ with basis $e_1,e_2,e_3$, put
$h=-e_1-e_2-e_3$, and let $S=\{e_1,e_2,e_3,h\}$.
A sequence $e_1^{a_1}e_2^{a_2}e_3^{a_3}h^b$ is a block
precisely when $a_i\equiv b\pmod4$ for each $i$. An atom with
any exponent at least four is one of the pure blocks
$e_1^4,e_2^4,e_3^4,h^4$. If all exponents are less than four,
they must be equal; the only atom of this kind is
$e_1e_2e_3h$. All atoms therefore have sequence length four, and
$S$ is half-factorial.

In $2G\cong C_2^3$, the elements $2e_1,2e_2,2e_3$ are a basis
and $2h=2e_1+2e_2+2e_3$. Thus $2S$ is exactly the configuration
of Example~\ref{ex:weak}. It is not half-factorial. Since any
subset of a half-factorial set is half-factorial, the set
$S\cup2S$ cannot be half-factorial either.
\end{proof}

There is a separate obstruction to extending prescribed generators
of a subgroup. The example for $p=2$ below appears in
\cite[p.~851, Remark~2]{GeroldingerGoebel}. The same construction
works for every prime, with an explicit atom witnessing failure.

\begin{proposition}\label{prop:extension}
Let $p$ be a prime and let
$G=C_p\oplus C_{p^2}=\langle a\rangle\oplus\langle e\rangle$,
where $\ord(a)=p$ and $\ord(e)=p^2$. Put
$b=a+pe$, $H=G[p]$, and $T=\{a,b\}$. Then $T$ is a
half-factorial generating set of the proper subgroup $H$, but no
half-factorial subset of $G$ containing $T$ generates $G$.
More precisely, for every $x\in G\setminus H$ there is an atom
of $\B(T\cup\{x\})$ with cross number $1+1/p$.
\end{proposition}

\begin{proof}
The elements $a,b$ are independent of order $p$ and generate
$H=\langle a,pe\rangle$, so $T$ is half-factorial. Let
$x\in G\setminus H$. Then $\ord(x)=p^2$ and
$px=jpe=j(b-a)$ for some $j\in\{1,\ldots,p-1\}$. The sequence
\begin{equation}\label{eq:finite-obstruction}
 W=x^p a^j b^{p-j}
\end{equation}
is a block. In a zero-sum subsequence the multiplicity of $x$
is divisible by $p$, by projection to $G/H$. If it is zero,
independence of $a,b$ forces the subsequence to be empty. If it
is $p$, the coefficients of $a,b$ must be $j,p-j$, respectively.
Thus $W$ is an atom. Its cross number is
$p/p^2+j/p+(p-j)/p=1+1/p$, so
Proposition~\ref{prop:cross} excludes half-factoriality as soon
as any $x\notin H$ is adjoined.
\end{proof}

By contrast, the union of a chain of \emph{nested} half-factorial
subsets is half-factorial: any two factorizations witnessing failure
would involve only finitely many letters, all belonging to one
member of the chain. This finite-character observation is standard
\cite[Lemma~3.1]{GaoGeroldinger}. It does not provide the nested
sets. Proposition~\ref{prop:extension} shows why compatibility at
successive stages requires its own argument.

The next example strengthens the warning following
Theorem~\ref{thm:weak-existence}. Failure is not confined to taking
too many generators: the weakly half-factorial set can be
irredundant. A generating set is \emph{irredundant} if none of its
proper subsets generates the group.

\begin{theorem}\label{thm:no-thinning}
Let $p$ be a prime, let $m\ge2$, and put $r=mp$. Let
\[
 E=\bigoplus_{i=1}^{r}\langle e_i\rangle,
 \qquad \ord(e_i)=p^2,
 \qquad G=E/\langle p(e_1+\cdots+e_r)\rangle.
\]
Write $t_i$ for the image of $e_i$ and let $T=\{t_1,\ldots,t_r\}$.
Then $T$ is an irredundant weakly half-factorial generating set of
$G$, but contains no half-factorial generating subset. More
precisely, every $t_i$ has order $p^2$, the atoms over $T$ are exactly
\begin{equation}\label{eq:thinning-atoms}
 t_1^{p^2},\ldots,t_r^{p^2},
 \qquad U=\prod_{i=1}^{r}t_i^p,
\end{equation}
and $\mathsf k(U)=m$. For every integer $q\ge1$,
\begin{equation}\label{eq:thinning-lengths}
 \LL_{\B(T)}(U^{pq})
   =\{pq+(r-p)j:0\le j\le q\}.
\end{equation}
Nevertheless, $G\cong C_p\oplus C_{p^2}^{\,r-1}$, so $G$ itself
has a half-factorial generating set.
\end{theorem}

\begin{proof}
The homomorphism on $E$ sending every $e_i$ to $1/p^2+\Z$
annihilates $p(e_1+\cdots+e_r)$, whose image is $r/p+\Z=m+\Z$.
It therefore descends to a character $\chi$ of $G$. In particular,
$\ord(t_i)=p^2$ and $\chi(t_i)=1/\ord(t_i)+\Z$. The $t_i$
are distinct: a relation $e_i-e_j$ cannot be a multiple of
$p(e_1+\cdots+e_r)$, as is already clear after reduction modulo
$pE$. Lemma~\ref{lem:character} proves weak half-factoriality.

The defining relation lies in $pE$, so
$G/pG\cong E/pE\cong C_p^{\,r}$, with the images of the $t_i$
forming a basis. Consequently no proper subset of $T$ generates
$G$.

A sequence $\prod_i t_i^{v_i}$ is a block precisely when, for one
$j\in\{0,\ldots,p-1\}$, all its multiplicities satisfy
$v_i\equiv jp\pmod{p^2}$. If any multiplicity is at least $p^2$,
an atom must be the corresponding pure block. Otherwise all
multiplicities equal $jp$. A nonempty such block is $U^j$, and
it is an atom exactly when $j=1$. This proves
\eqref{eq:thinning-atoms}. In particular,
$\mathsf k(U)=rp/p^2=m\ne1$, so $T$ is not half-factorial.
Together with irredundancy, this proves the assertion about all
its subsets.

A factorization of $U^{pq}$ uses some number $a$ of copies of
$U$ and numbers $b_i$ of the pure atoms. Comparing multiplicities
gives $pa+p^2b_i=p^2q$ for every $i$. Thus $a=pc$ and
$b_i=q-c$ with $0\le c\le q$. The factorization length is
$pc+r(q-c)$, giving \eqref{eq:thinning-lengths}. For example,
$U^p$ has lengths $p$ and $r=mp$.

Finally, $e_1+\cdots+e_r,e_2,\ldots,e_r$ is a basis of the
free $\Z/p^2\Z$-module $E$. In this basis the only additional
relation is that $p$ times its first member is zero. This gives
the asserted cyclic decomposition, whose standard generators
are half-factorial.
\end{proof}

This phenomenon persists among reduced uncountable unbounded
$p$-groups. For any uncountable cardinal $\kappa$,
take the direct sum of $\kappa$ copies of the group in
Theorem~\ref{thm:no-thinning}, together with
$\bigoplus_{n\ge1}C_{p^n}$. The union of the displayed weak
generating sets and the cyclic generators is irredundant and
weakly half-factorial, but has no half-factorial generating
subset. The group is reduced, uncountable, and unbounded.
It is nonetheless a direct sum of cyclic groups and belongs to
$\HF$. Thus the obstruction to thinning is not a counterexample
to the universal primary-group realization problem.

\subsection{The full character set in the torsion-completion}

One might instead keep all character-compatible elements, so that
no useful generator is discarded. The following elementary
construction shows why this also fails in groups with sufficiently
large socle. Its conclusion concerns the full set in
\eqref{eq:admissible-set}, not every generating subset of that set.

\begin{proposition}\label{prop:full-character-atoms}
Let $p$ be a prime, let $P$ be an abelian $p$-group, and let $\chi:P\to\Q/\Z$ be a
homomorphism with $\chi(P[p])\ne0$. If
$\dim_{\mathbb F_p}(P[p]\cap\ker\chi)$ is infinite, then, for
every integer $m\ge1$, the block monoid
$\B(W_\chi(P))$ has an atom of cross number $m$.
For each $m\ge2$ it has an element whose set of lengths contains
both $p$ and $mp$. In particular, the ratios of maximal to
minimal factorization lengths in this block monoid are unbounded.
\end{proposition}

\begin{proof}
Choose $c\in P[p]$ with $\chi(c)=1/p+\Z$. The block $c^p$
is an atom of cross number one. For $m\ge2$, put $r=mp-1$
and choose independent elements
$a_1,\ldots,a_r\in P[p]\cap\ker\chi$. Define
\[
 s_0=c-\sum_{i=1}^r a_i,\qquad
 s_i=c+a_i\quad(1\le i\le r).
\]
The elements $c,a_1,\ldots,a_r$ are independent over
$\mathbb F_p$. Hence the $s_i$ are distinct elements of order
$p$ with character value $1/p+\Z$, and all belong to
$W_\chi(P)$. The sequence $U=s_0s_1\cdots s_r$ has sum
$(r+1)c=mpc=0$.

If a subsequence of $U$ has multiplicities
$\varepsilon_i\in\{0,1\}$ and sum zero, its $a_i$-coordinate
forces $\varepsilon_i=\varepsilon_0$ for each $i$. Therefore
that subsequence is empty or is all of $U$. Thus $U$ is an atom
with $\mathsf k(U)=(r+1)/p=m$. Finally,
\[
 U^p=\prod_{i=0}^r s_i^p
\]
has a factorization of length $p$ into copies of $U$ and a
factorization of length $r+1=mp$ into pure atoms. Their length
ratio is $m$, which proves the last assertion.
\end{proof}

For the torsion-completion $P$ in
Proposition~\ref{prop:completion}, this applies to every character
for which $W_\chi(P)$ generates $P$. Indeed, such a character is
nonzero on $P[p]$: if $t\in W_\chi(P)$ has order $p^n$, then
$\chi(p^{n-1}t)=1/p+\Z$. Also,
$P[p]=\prod_{n\ge1}C_n[p]$ is infinite-dimensional, so the kernel
of its nonzero functional to the subgroup of order $p$ has
infinite dimension. Thus the full character set always has
atoms of arbitrarily large integral cross number in this test
case. This proves the failure of the full-set construction;
nonexistence of every half-factorial generating subset follows
instead from Theorem~\ref{thm:counterexample}.

The logical implications can now be summarized without confusing a
property of a specified set with existence of some suitable set.

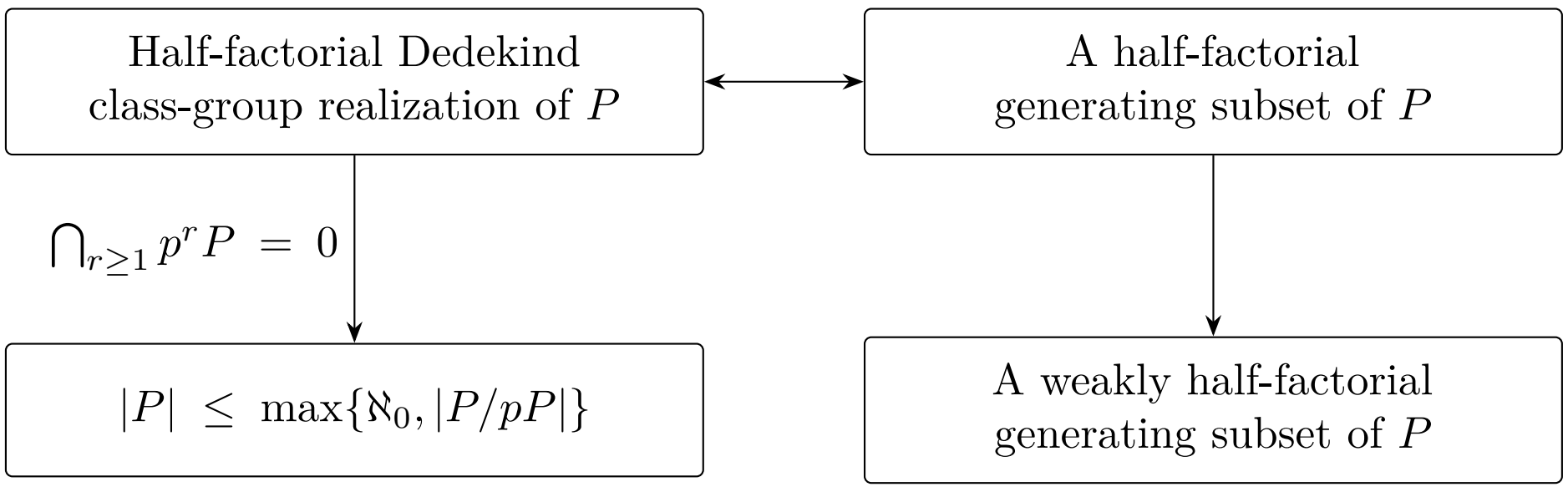
\begin{figure}[htbp]
\centering
\begin{tikzpicture}[>=Stealth,
 box/.style={draw,rounded corners=1.5pt,align=center,
             text width=5.1cm,minimum height=1.05cm,inner sep=6pt},
 every node/.style={font=\small}]
\node[box] (real) at (-3.4,0)
 {Half-factorial Dedekind\\class-group realization of $P$};
\node[box] (hf) at (3.4,0)
 {A half-factorial\\generating subset of $P$};
\node[box] (bound) at (-3.4,-2.6)
 {$|P|\le\max\{\aleph_0,|P/pP|\}$};
\node[box] (weak) at (3.4,-2.6)
 {A weakly half-factorial\\generating subset of $P$};
\draw[<->] (real) -- (hf);
\draw[->] (hf) -- (weak);
\draw[->] (real) -- node[left,align=right,text width=2.9cm]
 {$\bigcap_{r\ge1}p^rP=0$} (bound);
\end{tikzpicture}
\caption{Implications for an abelian $p$-group $P$. The equivalence is
classical; the left downward implication requires the displayed
separatedness hypothesis. Weakly half-factorial generators always
exist, but they do not imply half-factorial realizability. No
converse to the cardinality implication is asserted.}
\label{fig:implications}
\end{figure}

\section{Socle compatibility and repeated contraction}\label{sec:socle}

We finish with two structural restrictions on chosen generating
sets. At the prime two, containing a socle basis forces compatibility
with the entire finite-height filtration. Under independence modulo
$p$, one contraction preserves every set of lengths, but this
property need not persist under iteration. These results provide
finer local information than the global cardinality obstruction.

\subsection{A socle-filtration obstruction at the prime two}

For an abelian $2$-group $P$, its socle $P[2]$ is a vector space
over $\mathbb F_2$. A subset of $P[2]\setminus\{0\}$ is
half-factorial precisely when it is linearly independent. Indeed,
independence leaves only the pure atoms $e^2$. A dependent subset
has a finite minimally dependent subset of at least three distinct
elements, whose product is an atom of cross number at least $3/2$.

The next estimate bounds the socle support of a sequence by its
cross number. It applies to arbitrary sequences, not only to blocks.

\begin{lemma}\label{lem:socle-support}
Let $P$ be an abelian $2$-group, let $S\subseteq P\setminus\{0\}$
be half-factorial, and suppose that $S$ contains a vector-space
basis $E$ of $P[2]$. For $v\in P[2]$, let $d_E(v)$ be the number
of nonzero coordinates of $v$ in the basis $E$. If
$A\in\F(S)$ and $\sigma(A)\in P[2]$, then
\begin{equation}\label{eq:socle-support-bound}
 d_E(\sigma(A))\le 2\mathsf k(A).
\end{equation}
In particular, if $\mathsf k(A)=1/2$ and
$\sigma(A)\in P[2]\setminus\{0\}$, then $\sigma(A)\in E$.
\end{lemma}

\begin{proof}
Write $\sigma(A)=e_1+\cdots+e_d$ with distinct $e_i\in E$,
and put $C=e_1\cdots e_d$. Since each $e_i$ has order two,
$AC$ is a block. Mark the $d$ occurrences belonging to $C$ and
factor $AC$ into atoms. No atom contains two marked occurrences:
they are distinct letters of cross number $1/2$ each, so two of
them already have cross number one and do not sum to zero. Any
additional letter would give cross number greater than one,
contrary to Proposition~\ref{prop:cross}.

Consequently the factorization has at least $d$ atoms. Its length
is $\mathsf k(AC)=\mathsf k(A)+d/2$, again by
Proposition~\ref{prop:cross}. Thus $\mathsf k(A)+d/2\ge d$,
which proves \eqref{eq:socle-support-bound}. The final assertion
follows because a nonzero socle element has positive support size.
\end{proof}

For example, if $s\in S$ has order $2^r$, applying the lemma to
$s^{2^{r-1}}$ shows that $2^{r-1}s\in E$. The same observation
applies to a sequence $A$ with $\mathsf k(A)=2^{-r}$ whenever
$\sigma(A)$ has order $2^r$. This allows a carrying procedure
that respects the socle basis.

\begin{theorem}\label{thm:socle-filtration}
Let $P$ be an abelian $2$-group and let
$S\subseteq P\setminus\{0\}$ be a half-factorial generating
set. Suppose that $E\subseteq S$ is a vector-space basis of
$P[2]$. Put $V_n=(2^nP)\cap P[2]$ for $n\ge0$. Then
\begin{equation}\label{eq:socle-adapted}
 V_n=\operatorname{span}_{\mathbb F_2}(E\cap V_n)
 \qquad(n\ge0).
\end{equation}
Thus the same basis $E$ must be adapted to every finite-height
subspace of the socle.
\end{theorem}

\begin{proof}
\emph{Step 1: establish the invariants for a carrying procedure.}
We first record an elementary consequence of half-factoriality.
If a nonempty sequence $A$ has cross number $2^{-r}$, with
$r\ge1$, and $2^r\sigma(A)=0$, then $\sigma(A)$ has order
exactly $2^r$. Otherwise, for its order $d<2^r$, the nonempty
block $A^d$ would have cross number $d/2^r<1$, which is impossible.
When its order is $2^r$, Lemma~\ref{lem:socle-support}, applied
to $A^{2^{r-1}}$, gives $2^{r-1}\sigma(A)\in E$.

Fix $n\ge0$ and $v\in V_n\setminus\{0\}$. Choose $z\in P$
with $2^nz=v$; then $2^{n+1}z=0$. Since $S$ generates $P$ as
a monoid, choose a sequence $W\in\F(S)$ with $\sigma(W)=z$.
Initially partition $W$ into its individual letters. Each piece
$A$ of this partition satisfies, for some $r\ge1$,
\begin{equation}\label{eq:piece-invariant}
 \mathsf k(A)=2^{-r},\qquad
 \ord(\sigma(A))=2^r,\qquad
 2^{r-1}\sigma(A)\in E.
\end{equation}
The last condition follows from Lemma~\ref{lem:socle-support}.

\emph{Step 2: lower the largest order by pairing.}
Suppose that the largest exponent $r$ among the pieces exceeds
$n+1$. Multiplying the sum of all pieces by $2^{r-1}$ gives zero.
Pieces of smaller order disappear, while each piece of order
$2^r$ contributes its indicated basis element of $E$.
Independence of $E$ therefore shows that, for each $e\in E$,
an even number of the order-$2^r$ pieces have top multiple $e$.
Pair those pieces within each such group.

If $A$ and $A'$ are paired, then $AA'$ has cross number
$2^{-(r-1)}$ and its sum is annihilated by $2^{r-1}$, since the
two top multiples are equal. As $r-1\ge1$, the observation at
the start of the proof shows that its sum has order exactly
$2^{r-1}$; the same observation and the socle support lemma
establish the last condition of \eqref{eq:piece-invariant} for
the new piece. Replace every pair by its concatenation. The
largest exponent decreases. Repetition therefore yields a finite
partition of $W$ into pieces satisfying \eqref{eq:piece-invariant}
with $r\le n+1$.

\emph{Step 3: recover the required socle expansion.}
Finally multiply their sums by $2^n$. The pieces with $r\le n$
vanish. Every remaining piece contributes a member of $E$ that
belongs to $2^nP$, hence to $E\cap V_n$. Their sum is $v$.
This proves the nontrivial inclusion in \eqref{eq:socle-adapted};
the reverse inclusion is immediate.
\end{proof}

The theorem gives an obstruction when the socle filtration has
small successive quotients but a large underlying vector space.
Indeed, if $\bigcap_{n\ge0}V_n=0$, then every $e\in E$ lies
in exactly one set $E\cap(V_n\setminus V_{n+1})$. The images
of this set form a basis of $V_n/V_{n+1}$ by
\eqref{eq:socle-adapted}. In particular, countable-dimensional
successive quotients would force $E$ to be countable.

\begin{corollary}\label{cor:no-nested-completion}
Let
\[
 P=t\!\left(\prod_{j\ge1}\Z/2^j\Z\right),
\]
where $t(-)$ denotes the torsion subgroup. There is no
half-factorial generating subset $S\subseteq P\setminus\{0\}$
for which $S\cap P[2]$ generates $P[2]$.
Consequently, there are no subsets
$T_n\subseteq P[2^n]\setminus\{0\}$, $n\ge1$, such that
all $T_n$ are half-factorial, all $T_n$ generate $P[2^n]$, and
$T_n\subseteq T_{n+1}$ for every $n$.
\end{corollary}

\begin{proof}
Identify $P[2]$ with $\prod_{j\ge1}\mathbb F_2$ by writing its
$j$th coordinate as $\varepsilon_j2^{j-1}+2^j\Z$.
For $V_n=(2^nP)\cap P[2]$ one has
\begin{equation}\label{eq:completion-filtration}
 V_n=\{(\varepsilon_j)_j:\varepsilon_1=\cdots=
 \varepsilon_n=0\}.
\end{equation}
The forward inclusion is immediate. For the reverse inclusion,
choose the $j$th coordinate of a $2^n$th root to be
$\varepsilon_j2^{j-n-1}+2^j\Z$ when $j>n$, and zero when
$j\le n$. This root is annihilated by $2^{n+1}$ and therefore
belongs to $P$.

It follows that $\bigcap_nV_n=0$ and
$\dim_{\mathbb F_2}(V_n/V_{n+1})=1$ for every $n\ge0$.
On the other hand, $P[2]=\prod_{j\ge1}\mathbb F_2$ is
uncountable and consequently has no countable vector-space basis.

Suppose that $S$ as in the first assertion existed. The subset
$E=S\cap P[2]$ would be half-factorial and would generate the
socle. By the observation at the start of this subsection it
would be a basis of $P[2]$. Theorem~\ref{thm:socle-filtration}
and \eqref{eq:completion-filtration} imply that
$E\cap(V_n\setminus V_{n+1})$ has exactly one member for each
$n$. Their union is $E$, since $\bigcap_nV_n=0$. This makes
$E$ countable, a contradiction.

For the last assertion, the union $S=\bigcup_nT_n$ would be
half-factorial by finite character and would generate
$P=\bigcup_nP[2^n]$. Its intersection with $P[2]$ contains
$T_1$ and hence generates the socle, contradicting the first
assertion.
\end{proof}

The conclusion is stronger than Example~\ref{ex:incompatible-layers}:
that example exhibits unsuccessful choices, whereas
Corollary~\ref{cor:no-nested-completion} excludes every nested
choice in the indicated group. The socle argument by itself proves only a conditional obstruction;
the unconditional nonrealization result is
Theorem~\ref{thm:counterexample}. A global generating set need
not generate the socle using its order-two members.
The next elementary example makes the distinction explicit.

\begin{example}\label{ex:socle-not-generated}
Let $G=C_2\oplus C_4=\langle a\rangle\oplus\langle z\rangle$,
where $\ord(a)=2$ and $\ord(z)=4$, and put $S=\{z,z+a\}$.
This is a generating set and both its elements have order four.
Its block atoms are exactly
\[
 z^4,\qquad (z+a)^4,\qquad z^2(z+a)^2.
\]
Indeed, a sequence $z^u(z+a)^v$ is a block precisely when
$v$ is even and $u+v$ is divisible by four, from which the list
follows. Every listed atom has cross number one. Thus $S$ is
half-factorial, although $S\cap G[2]=\varnothing$ and
$G[2]\ne0$. Therefore the socle-generation hypothesis in
Corollary~\ref{cor:no-nested-completion} cannot be inferred merely
from half-factoriality and generation of the whole group.
\end{example}

\subsection{Contraction under independence modulo a prime}

A different possible approach is to contract generators repeatedly
by multiplication by $p$. The following criterion explains when
one contraction preserves all factorization lengths. It does not
assert preservation for an arbitrary generating set.

\begin{proposition}\label{prop:independent-contraction}
Let $p$ be a prime, let $P$ be an abelian $p$-group, and let
$S\subseteq P\setminus\{0\}$ generate $P$. Suppose that the
images of the members of $S$ are linearly independent in $P/pP$.
Every $B\in\B(S)$ then has all its multiplicities divisible
by $p$, and the map
\[
 \theta:\B(S)\longrightarrow\B(pS),\qquad
 \prod_{s\in S}s^{p a_s}\longmapsto
 \prod_{s\in S}(ps)^{a_s}
\]
is a surjective homomorphism preserving every set of factorization
lengths. Here coincident letters on the right are combined, and
$0$ is allowed in $pS$. In particular, $S$ is half-factorial
if and only if $pS$ is half-factorial.
\end{proposition}

\begin{proof}
Reducing a zero-sum relation modulo $pP$ and using the stated
independence shows that each multiplicity is divisible by $p$.
The displayed substitution consequently defines a homomorphism
into $\B(pS)$. Given a block over $pS$, choose a preimage in $S$
for each occurrence and replace that occurrence by $p$ copies of
its preimage. The resulting sequence is a block over $S$, proving
surjectivity.

For a given $B\in\B(S)$, group its occurrences of each $s$
into packets of $p$ identical letters, retaining their labels in
$\theta(B)$. Any factorization of $\theta(B)$ can be lifted by
distributing those labelled packets among the factors. Conversely,
if the image of an atom factored, this lifting procedure would
factor the original atom. Thus atoms map to atoms, and every
factorization of the image lifts to a factorization with the same
length. This proves the asserted equality of length sets. Adding
or deleting the letter zero does not affect half-factoriality.
\end{proof}

The independence needed for the next contraction can disappear.
In fact, arbitrarily many successful contractions can precede the
first failure, even when the original generating set is
irredundant and the group is finite.

\begin{theorem}\label{thm:delayed-contraction}
Let $p$ be a prime and $r\ge1$. Put $k=p^2$, $N=r+2$, and
\[
 G=\left(\bigoplus_{i=1}^{k}C_{p^N}e_i\right)
 \Big/\left\langle p^r(e_1+\cdots+e_k)\right\rangle.
\]
Let $t_i$ be the image of $e_i$ and set $S=\{t_1,\ldots,t_k\}$.
Then $S$ is an irredundant half-factorial generating set of $G$,
and $p^jS$ is half-factorial for every $0\le j\le r$, but
$p^{r+1}S$ is not half-factorial. More precisely,
$\prod_{i=1}^k(p^{r+1}t_i)$ is an atom of cross number $p$.
Moreover, $G\cong C_{p^r}\oplus C_{p^{r+2}}^{\,p^2-1}$.
\end{theorem}

\begin{proof}
Each $t_i$ has order $p^N$. Indeed, a multiple of $e_i$ can
belong to the displayed relation subgroup only if all the other
coordinates of that relation vanish, in which case its $i$th
coordinate also vanishes. The relation vector belongs to $p$
times the original direct sum, so the images of the $t_i$ form a
basis of $G/pG$. In particular, $S$ is irredundant.

For $0\le j\le r+1$ put $x_i=p^jt_i$,
$q=p^{N-j}$, and $d=p^{\max\{r-j,0\}}$. The $x_i$ are
nonzero and distinct, and all have order $q$. Direct inspection
of the relation subgroup shows that $\prod_i x_i^{a_i}$ is a
block precisely when all $a_i$ have the same residue modulo $q$
and that common residue is divisible by $d$. Therefore its atoms
are exactly the pure blocks $x_i^q$ and the diagonal block
\[
 D_j=(x_1\cdots x_k)^d.
\]
To verify completeness of this list, an atom with an exponent at
least $q$ is a pure block. Otherwise all exponents are the same
integer $\ell d$ in $[0,q-1]$, and a nonempty such block is a
power of $D_j$.

If $j\le r$, then $\mathsf k(D_j)=kd/q=1$, so the
cross-number criterion proves half-factoriality. If $j=r+1$,
then $q=p$, $d=1$, and $\mathsf k(D_j)=k/p=p>1$.
For an explicit factorization witness,
\[
 D_{r+1}^{\,p}=\prod_{i=1}^{p^2}x_i^p
\]
has lengths $p$ and $p^2$. Finally, replace $e_k$ in the standard
basis of $\bigoplus_iC_{p^N}e_i$ by $e_1+\cdots+e_k$.
The relation then affects only that cyclic summand, reducing its
order to $p^r$ and giving the stated group decomposition.
\end{proof}

Thus one-step preservation under the hypotheses of
Proposition~\ref{prop:independent-contraction} does not justify
an unrestricted induction. In particular, successful testing
through any prescribed finite number of contractions supplies no
uniform guarantee for the next contraction. This obstruction, like
the socle obstruction above, delimits a construction method rather
than producing a group with no half-factorial generating set.

\subsection{Concluding remarks}

The separation theorem answers the universal Dedekind realization
problem negatively by controlling arbitrary global generators.
For each prime $p$, the torsion-completion $P_p$ has countable
finite power quotients but is uncountable. Half-factorial generators
would be a countable union of sets injecting into these quotients,
which is impossible. The construction of pure subgroups shows that
the obstruction already occurs at the least uncountable cardinal.

Several positive facts remain compatible with this negative result.
Every torsion abelian group has weakly half-factorial generators;
every bounded subgroup has half-factorial generators; and the
completion is a pure extension with both kernel and quotient in
$\HF$. None of these statements ensures a single global
half-factorial generating set. The failure of closure under pure
extensions and countable increasing unions makes the distinction
an algebraic one, not merely a failure of a particular proof method.

The cardinal bound is a necessary condition for $p$-adically
separated groups. Outside that hypothesis the proof does not apply,
since some nonzero top multiples may belong to every $p^rP$.
Even within the separated class, the bound alone is not proved
sufficient. A structural characterization of the groups belonging
to $\HF$ therefore remains a separate problem. The splitting,
character, socle, and contraction results offer restrictions and
constructions for that classification without being needed for the
negative answer established in Section~\ref{sec:cardinality}.

\bigskip

\paragraph{\bf Competing interests:} The authors declare that they have no competing interests.

\medskip

\paragraph{\bf Availability of data and material:}
No datasets were generated or analysed during the current study.

\medskip

\paragraph{\bf Use of large language models, AI, and machine learning tools.} 
During the preparation of this work the authors
used ChatGPT, an AI assistant developed by OpenAI, to explore various mathematical concepts related to the topics of this paper and to verify the English grammar. After using this tool/service,
the author reviewed and edited the content as needed and take full responsibility
for the content of the published article.

\medskip

\paragraph{\bf Funding}
H. Kim was supported by the Basic Science Research Program through
the National Research Foundation of Korea (NRF), funded by the Ministry
of Education (grant no.~RS-2026-25589739).

\end{document}